\documentclass{article}
\usepackage{amsmath, latexsym, amsfonts, amssymb, amsthm, amscd}

\newtheorem {lemme} {Lemma} [section]
\newtheorem {theoreme} {Theorem} [section]
\newtheorem {proposition} {Proposition} [section]

\newtheorem {remarque} {Remark} [section]

\newcommand{\tr}{{\rm tr}}
\newcommand{\Tr}{{\rm Tr}}
\newcommand{\E}{\mathbb {E}}

\newcommand{\R}{\mathbb {R}}
\newcommand{\C}{\mathbb {C}}

\newcommand{\1}{1\!\!{\sf I}}

\numberwithin{equation}{section}

\newcommand{\vers}{\mathop{\longrightarrow }}
\title{Additive/multiplicative free subordination property and limiting
eigenvectors of spiked additive deformations of Wigner matrices and spiked
sample covariance matrices
\footnote{This work was partially supported by the {\emph Agence Nationale de la
Recherche} grant ANR-08-BLAN-0311-03.}}

\author{M. Capitaine\thanks{CNRS, Institut de Math\'ematiques de Toulouse, 
Equipe de Statistique et Probabilit\'es,  F-31062 Toulouse Cedex 09. 
E-mail: mireille.capitaine@math.univ-toulouse.fr }}

\date{}
\begin{document}
\maketitle
\begin{abstract}
When  some eigenvalues of  a spiked additive deformation of a   Wigner matrix or a spiked  multiplicative deformation of a Wishart matrix  separate from the bulk, we  study  how the corresponding eigenvectors  project  onto those of the perturbation. 
We point out that  the   subordination function relative to the free (additive or multiplicative) convolution plays an important part in the asymptotic behavior.
\end{abstract}

\section{Introduction}\label{intro}
This paper lies in the lineage of recent works studying the
influence of some perturbations on the asymptotic spectrum of
classical random matrix models.
Such questions come from Statistics
(cf. \cite{John}) and appeared in the framework of empirical
covariance matrices.
In the pioneering work \cite{BBP05}, J. Baik, G. Ben Arous and S.
P\'ech\'e  dealt with random sample covariance matrices
$(S_N)_N$ defined by
\begin{eqnarray}{\label{spike}}
S_N=\frac{1}{p} Y_N Y_N^* \mbox{~~with~~} Y= \Sigma_N^{\frac{1}{2}} B_N
\end{eqnarray}
where $B_N$ is a $N \times p$ complex matrix such that the entries $(B_N)_{ij}$ are i.i.d centered standard Gaussian and ${\Sigma}_N$
is a deterministic positive $N\times N$ matrix  having all but finitely
many eigenvalues equal to one. This model can be seen as a multiplicative perturbation of the so-called white Wishart matrix for which ${\Sigma}_N=I_N$. Besides, the size of the samples $N$
and the size of the population $p=p_N$ are assumed of the same order
(as $N \to \infty$, $N/p \rightarrow c>0$). The global limiting behavior
of the spectrum of $S_N$ is not affected by such a matrix $\Sigma_N$.
Thus,  the limiting spectral measure is the well-known
Marchenko-Pastur law (\cite{MP67})
defined by
\begin{equation} \label{MP} \mu_{\mbox{\tiny{MP}},{c}}(dx)=\max\{1-\frac{1}{c},0\}\delta_0
+f(x) \1_{[(1-\sqrt{c})^2;
(1+\sqrt{c})^2]}(x)dx\end{equation}
with $$ f(x)=\frac{\sqrt{\left(x-(1-\sqrt{c})^2\right)
\left((1+\sqrt{c})^2-x\right)}}{2\pi c x}.$$
When ${\Sigma}_N=I_N$, the largest eigenvalue of $S_N$ converges towards the right hand point of the support of the Marchenko-Pastur law (see \cite{Ge,BaiSilYi,YBK}). When ${\Sigma}_N\neq I_N$,
in \cite{BBP05} the authors pointed out a striking phase
transition phenomenon for the asymptotic behaviour of the largest eigenvalue of $S_N$ (at the convergence and fluctuations levels)
according to the value of the largest eigenvalue(s) of $\Sigma_N$.
They showed in particular that when
the largest eigenvalue of ${\Sigma}_N$ is far from one, the largest eigenvalue of $S_N$ converges   outside the limiting
Marchenko-Pastur support. In  \cite{BaikSil06},
 J. Baik and
J. Silverstein extended the result of \cite{BBP05} on  the convergence of the  extremal eigenvalues of complex or real non
necessarily Gaussian matrices $S_N$ under  finite   four
moments assumptions on the distribution of the entries of $B_N$. \\ 

When $S_N$ is still defined by 
(\ref{spike}),
but now the limiting spectral distribution of  $\Sigma_N$ is some compactly supported measure $\nu$ on $[0;+ \infty[$, under  finite  second
moments assumptions on the distribution of the entries of $B_N$,  the  spectral distribution of $S_N$ converges almost surely towards 
a probability measure $\mu_{\tiny{LSD}}$ which only depends on $c$ and $\nu$; denoting by 
$g_{\mu_{\tiny{LSD}}}(z)= \int \frac{1}{z-x} d\mu_{\tiny{LSD}}(x)$ the  Stieltjes transform of  $\mu_{\tiny{LSD}}$,  
 for $z \in \mathbb{C}^+$,  $g_{\mu_{\tiny{LSD}}}(z)$ is the unique solution $Z$ in $ \{ Z \in \mathbb{C}, -\frac{(1-c)}{z} -cZ \in \mathbb{C}^+ \} $  of the equation  
\begin{equation}\label{MPeq} Z=\int \frac{1}{z- t(1-c +czZ)}d\nu(t) \end{equation}
 (see \cite{MP67,BaiSil95,GS,KY,SSt,Wa,Yin86}).
Very recently R. Rao and J. Silverstein \cite{RaoSil09} and 
 Z.~D. Bai and J. Yao \cite{BaiYao08b} dealt with such a model assuming  moreover that  $\Sigma_N$
 has a finite number of eigenvalues fixed outside the support of  $\nu$ called spikes (or converging outside the support of  $\nu$ in \cite{RaoSil09}),  whereas the distance between the other eigenvalues of  $\Sigma_N$ and the support of $\nu$ uniformly goes to zero. (Note that the assumptions in \cite{RaoSil09} are a bit more general).
Under  finite   four
moments assumptions on the distribution of the entries of $B_N$, the authors characterized the spikes of $\Sigma_N$ that  generate jumps of eigenvalues of $S_N$ 
 and described the corresponding limiting points 
outside the support of the limiting spectral distribution $\mu_{\tiny{LSD}}$ of $S_N$.\\

Several authors  considered an additive analogue  of the above setting that is, the
influence on the asymptotic spectrum of
the addition of  some Hermitian deterministic perturbation $A_N$ to the rescaled so-called Hermitian Wigner $N\times N$ matrix $W_N$.

\noindent Recall that, according to Wigner's work \cite{Wigner55,Wigner58} 
and further results of different authors (see \cite{Bai99} for a review), 
provided the common distribution $\mu $ of the entries is centered
with variance $\sigma ^2$, the large $N$-limiting spectral distribution 
of the rescaled complex Wigner matrix $X_N=\frac{1}{\sqrt{N}} W_N$ is 
the semicircle distribution $\mu _{\sigma }$ whose density is given by
\begin{equation}\label{scl}
\frac{d\mu _{\sigma }}{dx}(x)= \frac{1}{2 \pi \sigma ^2} \sqrt{4\sigma ^2- x^2} 
\, 1 \hspace{-.20cm}1_{[-2\sigma , 2\sigma ]}(x).
\end{equation}
Moreover, if the fourth moment of the measure $\mu $ is finite, 
the largest (resp. smallest) eigenvalue of $X_N$ 
converges almost surely towards the right (resp. left) endpoint $2\sigma $ 
(resp. $-2\sigma $) of the semicircular support 
(cf. \cite{BaiYin88} or Theorem 2.12 in \cite{Bai99}). 

\noindent Let $A_N$ be a deterministic Hermitian matrix such that  the spectral measure of $A_N$ weakly converges to some probability measure $\nu$ 
and $\left\| A_N\right\| $ is uniformly bounded in $N$.
When $N$ becomes large, free probability provides us a good understanding 
of the global behaviour of the spectrum of $M_N=X_N+A_N$ where $X_N$ is a rescaled complex Wigner matrix. 
Indeed,
the spectral distribution of $M_N$ weakly converges to the free convolution $\mu _\sigma \boxplus \nu $ 
almost surely and in expectation 
(cf \cite{AGZ09,MinSpe10} and \cite{Voiculescu91,Dykema93} for pioneering works). 
We refer the reader to \cite{VDN92} for an introduction to free probability theory. 

\noindent Dealing with small rank perturbation of a G.U.E matrix $W_N^G$,
S. P\'ech\'e pointed out an analogue phase transition phenomenon as in the sample covariance setting
for the convergence and the fluctuations of the largest eigenvalue of $M_N^G=W_N^G/{\sqrt{N}}+A_N$ 
with respect to the largest eigenvalue $\theta $ (independent of $N$) of $A_N$ \cite{Peche06}. 
These investigations imply that, if $\theta $ is far enough from zero ($\theta > \sigma $), 
then the largest eigenvalue of $M_N^G$ jumps above the support $[-2\sigma , 2\sigma ]$ 
of the limiting spectral measure and converges (in probability) 
towards $\rho _\theta =\theta + \frac{\sigma ^2}{\theta }$. 
Note that Z. F\"uredi and J. Koml$\acute{\text{o}}$s already exhibited such a phenomenon
in \cite{FurKom81} dealing with non-centered symmetric matrices.\\
In \cite{F�P�07}, D. F\'eral and S. P\'ech\'e proved that the results of \cite{Peche06} 
still hold for a non-necessarily Gaussian Wigner Hermitian matrix $W_N$ 
with sub-Gaussian moments and in the particular case of a rank one perturbation matrix $A_N$ 
whose entries are all $\frac{\theta }{N}$ for some real number $\theta $. 
In \cite{CDF09}, the authors considered a deterministic Hermitian matrix $A_N$ 
of arbitrary fixed finite rank $r$ 
and built from a family of $J$ fixed non-null real numbers 
$\theta _1 > \cdots > \theta _J$ independent of $N$ 
and such that each $\theta _j$ is an eigenvalue of $A_N$ of fixed multiplicity $k_j$ 
(with $\sum_{j=1}^J k_j=r$). 
They dealt with general Wigner matrices associated to some symmetric measure satisfying a Poincar\'e inequality. 
They proved that eigenvalues of $A_N$ with absolute value strictly greater than $\sigma $
generate some eigenvalues of $M_N$ which converge 
to some limiting points outside the support of $\mu _\sigma $. 
In \cite{CDFF10}, the authors investigated the asymptotic behavior of the eigenvalues of generalized spiked perturbations of Wigner matrices associated to some symmetric measure satisfying a Poincar\'e inequality.
In this paper, the perturbation matrix $A_N$ is a deterministic Hermitian matrix whose spectral measure 
converges to some probability measure $\nu $ with compact support and such that $A_N$ has a fixed number of fixed eigenvalues (spikes) 
outside the support of $\nu $, whereas the distance between the other eigenvalues 
and the support of $\nu $ uniformly goes to zero as $N$ goes to infinity. 
It is established that only a particular subset of the spikes will generate 
some eigenvalues of $M_N$ which will converge to some limiting points 
outside the support of the limiting spectral measure. The phenomenon is completely analogous to the one described in \cite{RaoSil09} and 
 \cite{BaiYao08b} in the
sample covariance setting.\\

Now, one can wonder in the spiked deformed Wigner matrix setting as well as in the spiked sample covariance matrix setting, when some eigenvalues separate from the bulk, how the corresponding eigenvectors of the deformed model project onto those of the perturbation.
There are some results concerning finite rank perturbations: \cite{Paul} in the real Gaussian sample covariance matrix  setting, and \cite{BGRao09} dealing with finite rank additive or multiplicative perturbations of unitarily invariant matrices.
For a general perturbation, up to our knowledge nothing has been done concerning eigenvectors of deformed Wigner matrices.
Dealing with sample covariance matrices, S. P\'ech\'e and O. Ledoit \cite{PL} introduced a tool to study the average behaviour of the eigenvectors but it seems that this did not allow 
them  to focus on  the eigenvectors associated with  the eigenvalues that separate from the bulk.

As already said, the limiting spectral distribution of the deformed Wigner model is described by the free convolution
of the respective  limiting spectral distributions. Moreover,
the authors  explained in \cite{CDFF10} that the phenomenon of the eigenvalues separating from the bulk can be fully described in terms of free probability 
involving the subordination function related 
to the free additive convolution of a semicircular distribution with the limiting spectral distribution of the perturbation.
Actually, as we will show below, the analogue results in the sample covariance matrix  setting can also be
described in terms of free probability 
involving the subordination function related 
to the free multiplicative convolution of a Marchenko-Pastur distribution with the limiting spectral distribution of the perturbation.
Moreover, as already noticed by P. Biane in \cite{PB}, free probability again has something to tell us about eigenvectors of deformed matricial models. Indeed, in this paper,  we are going to describe in the deformed Wigner matrix setting as well as in 
 the sample covariance matrix one, how the eigenvectors of the deformed model associated to  the eigenvalues that separate from the bulk project onto those associated to the spikes of the perturbation, pointing out that  the  subordination functions relative to the free additive or  multiplicative convolution play an important part in this asymptotic behavior. Note that the proof is exactly the same in the additive and the multiplicative cases.\\

In the sample covariance matrix model as well as in the deformed Wigner model, the convergence of the eigenvalues that separate from the bulk is deduced from a striking exact separation phenomenon, roughly stating that to each gap in the spectrum of the deformed model there corresponds a gap in the spectrum of the perturbation, these gaps splitting in exactly the same way the corresponding spectrum. For general deformed models, that is, dealing with other matrices than Wigner matrices in the additive case or other matrices than white Wishart matrices in the multiplicative case, such an exact separation phenomenon is not expected in full generality. Nevertheless,
we express all the results in this paper in terms of the free additive respectively multiplicative subordination functions since we conjecture that, 
for other deformed models than deformed Wigner matrices and sample covariance matrices,
the limiting values of the eigenvalues that separate from the bulk as well as the limiting values of the orthogonal projection of the corresponding eigenvectors  onto those associated to the spikes of the perturbation will be given by the same quantities
 provided one deals with the corresponding subordination functions relative to the limiting spectral distribution of the non-deformed model. By the way, note that one can check that the results of F. Benaych-Georges and R. N. Rao in \cite{BGRao09}, concerning finite rank multiplicative or additive perturbation of a unitarily invariant matrix,
 about the convergence of the extremal eigenvalues and of the projection of the corresponding eigenvectors onto those of the perturbation can be rewritten in terms of subordination functions as conjectured.
 \\
 
 The paper is organized as follows. In Section 2, we introduce the additive and multiplicative deformed models we consider in this paper; we also introduce some basic notations that will be used throughout the paper.
 Section 3 is devoted to definitions and results concerning free convolutions and subordination functions, some of them being necessary to state our main result Theorem \ref{cvev1p1} in Section 4. Note that we present a common formulation for the additive and multiplicative deformed models
  and a common proof in Section 5 and Section 6, postponing in Section 7 the technical results that need a specific study for each model. Finally, an Appendix gathers several tools that will be used in the paper.

\section{Models and Notations}
Let $\mu$ be a probability measure with variance $\sigma ^2$ which  satisfies a Poincar\'e inequality with constant $C_{PI}$ (the definition of 
such an inequality is recalled in the Appendix). Note that this condition implies that $\mu$ has moments of any order (see Corollary 3.2 and Proposition 1.10 in \cite{Ledoux01}).
In this paper, we will deform  the 
following classical matricial models.
\begin{itemize}
\item Normalized Wigner matrices $X^W_N=\frac{1}{\sqrt{N}} W_N $ \\
such that
 $W_N$ is a $N \times N$ Wigner Hermitian matrix associated to 
the distribution $\mu $: \\
$(W_N)_{ii}$, $\sqrt{2} \Re ((W_N)_{ij})_{i < j}$, $\sqrt{2}
\Im ((W_N)_{ij})_{i<j}$
are  i.i.d., with distribution $\mu$. 
\item Sample covariance matrices $X^S_N=\frac{1}{p}B_N B_N^*$\\such that 
$B_N$ is a $N \times p$  matrix such that 
 $\sqrt{2} \Re ((B_N)_{ij})_{1\leq i\leq N, 1\leq j\leq p}$, $\sqrt{2}
\Im ((B_N)_{ij})_{{1\leq i\leq N, 1\leq j\leq p}}$
are  i.i.d., with distribution $\mu $. 
We assume that $\frac{N}{p} \rightarrow c>0$ when $N$ goes to infinity.
\end{itemize}
{\bf For Wigner matrices, we will assume moreover that $\mu$ is symmetric} since we will use results of \cite{CDFF10}
where this assumption is needed.\\

\noindent We will deform  these models by respectively addition and multiplication by a deterministic Hermitian perturbation matrix $A_N$; {\bf in the multiplicative perturbation case, $A_N$  will be  assumed to be nonnegative definite}.
 In both cases, we assume that:\\
 {\bf Assumption A}:\\
 The eigenvalues $\gamma_i=\gamma_i(N)$ of $A_N$
are such that the spectral measure $\mu _{A_N} := \frac{1}{N} \sum_{i=1}^N \delta _{\gamma _i}$ weakly
converges to some probability measure $\nu $ with compact support. 
We assume that there exists a fixed integer $r\geq 0$ (independent of $N$) 
such that $A_N$ has $N-r$ eigenvalues $\beta _j(N)$ satisfying 

 \begin{equation}\label{univconv}
 \max _{1\leq j\leq N-r} {\rm dist}(\beta _j(N),{\rm supp}(\nu ))\vers _{N \rightarrow \infty } 0,\end{equation}
where ${\rm supp}(\nu )$ denotes the support of $\nu $. \\
We also assume that there are $J$ fixed real numbers $\theta _1 > \ldots > \theta _J$ 
independent of $N$ which are outside the support of $\nu $ and such that each $\theta _j$ 
is an eigenvalue of $A_N$ with a fixed multiplicity $k_j$ (with $\sum_{j=1}^J k_j=r$). 
The $\theta_j$'s will be called {\bf the spikes or the spiked eigenvalues} of $A_N$.  The set of the spikes of $A_N$ will be denoted by
$\Theta$:
$$\Theta:=\{\theta_1; \ldots, \theta_J\}.$$ \\
{\bf In the sample covariance matrix setting we assume $\theta _J>0$.}\\

\noindent We will consider simultaneously the two deformed models:
$$M_N^W= X_N^W+A_N=\frac{1}{\sqrt{N}} W_N +A_N,$$
$$M_N^S=A_N^{\frac{1}{2}}X_N^S A_N^{\frac{1}{2}}=\frac{1}{p}A_N^{\frac{1}{2}}B_N B_N^*A_N^{\frac{1}{2}}.$$
{\bf When the approaches are the same for the two models we adopt the notation $M_N$ standing for both $M_N^W$
and $M_N^S$. When the studies are specific to one of the two models, we will use the superscripts.}\\

\noindent  Actually, we  assume without loss of generality in the sample covariance  setting that the variance of $\mu$ is 1
since it corresponds to considering the rescaled matrix $\frac{M^S_N}{\sigma^2}$.\\


\noindent Throughout this paper, we will use the following notations.
\begin{itemize}
\item[-] $\mathbb{C}^+$ will denote the complex upper half plane $\{z \in \mathbb{C}, \, \Im z >0\}$. Similarly,
$\mathbb{C}^-$ will stand for
$\{z \in \mathbb{C}, \, \Im z <0\}$.
\item[-] For a function $f$ differentiable in some neighborhood of a point $x$ in $\mathbb{R}$, we will denote by $f^{'}(x)$
the derivative of $f$.
\item[-] For a vector subspace ${\cal V}$ of $\mathbb{C}^N$, we will denote by ${\cal V}^{\bot}$ its orthogonal supplementary 
subspace
and by 
 $P_{\cal V}$ the orthogonal projection onto 
${\cal V}$.
\item[-]  $\langle ~, ~\rangle $ will denote the Hermitian inner product on $\mathbb{C}^N$ defined by $\langle a,b\rangle =b^*a$ for any $a,b$ in $\mathbb{C}^N$.
\item[-]  $\Vert ~ ~\Vert_2 $ will denote the Euclidean norm on $\mathbb{C}^N$.
\item[-] We will denote by $M_{m\times q}(\mathbb{C})$ the set of $ m\times q$ matrices with complex entries.
$\Vert  ~ \Vert $ will denote the operator norm
and  for any matrix $M$, $\Vert M \Vert_2= \{Tr(MM^*)\}^{\frac{1}{2}}$.
\item[-]  For any matrix $M$ in  $M_{N\times N}(\mathbb{C})$, we will denote its kernel by Ker($M$).
\item[-] $E_{ij}$ in $M_{m\times q}(\mathbb{C})$ stands for the matrix such that $(E_{ij})_{kl}=\delta_{ik}\delta_{jl}$.
\item[-] For any $N\times N$ Hermitian matrix $M$, we will denote by 
$$\lambda_1(M) \geq \ldots \geq \lambda_N(M)$$
\noindent its ordered eigenvalues.
\item[-] For a probability measure $\tau $ on $\R$, 
we denote by $\mbox{supp}(\tau)$ its support and by $^c \mbox{supp}(\tau) $ its complement in $\R$.
\item[-] For a probability measure $\tau $ on $\R$, 
we denote by $g_\tau $ its Stieltjes transform defined for $z \in \C\setminus \R$ by 
$$g_\tau (z) = \int_\R \frac{d\tau (x)}{z-x}.$$
\item[-] $G_N$ denotes the resolvent of $M_N$ and 
$g_N$ the mean of the Stieltjes transform of the spectral measure of $M_N$, that is, 
$$g_N(z) = \E(\tr_N G_N(z)), \, z \in \C\setminus \R,$$
where $\tr_N$ is the normalized trace: $\tr_N =\frac{1}{N} \Tr$.\\
When it is necessary to distinguish the deformed Wigner matrix setting and the sample covariance matrix one, we will specify 
the  resolvent or the Stieltjes transform  by using  the corresponding superscript as follows:
$G_N^W, g_N^W$ and $G_N^S, g_N^S$.
\item[-] $C, K $ denote nonnegative constants which may vary from line to line.
\end{itemize}
As already mentioned in the introduction, the assumptions on $W_N$ and $A_N$ ensure that  they are asymptotically free,
and then
the spectral distribution of $M_N^W$ weakly converges to the free convolution $\mu _\sigma \boxplus \nu $ 
almost surely and in expectation 
(cf \cite{AGZ09,MinSpe10} and \cite{Voiculescu91,Dykema93} for pioneering works). 

\noindent Concerning the sample covariance matrix model $M_N^S$, as already noticed in the introduction, its limiting spectral measure  only depends on $c$ and $\nu$.
Note that when the entries of $B_N$ are Gaussian (that is, if $\mu$ is Gaussian) we can assume that $A_N$ is diagonal by the invariance under unitary conjugation of the distribution of $X^S_N$. Then, since according to Corollary 4.3.8
in \cite{HP},
$X^S_N$ and $A_N$ are asymptotically free, we can conclude that {\bf the limiting spectral distribution of $M_N^S$ is actually  the free multiplicative convolution of the limiting spectral measure of $X_N^S$, that is, $ \mu_{\mbox{\tiny{MP}},{c}}$, with $\nu$, denoted by $ \mu_{\mbox{\tiny{MP}},{c}}\boxtimes \nu$.}

Thus, free additive and multiplicative convolutions provide a good understanding of the limiting global behaviour of the spectrum of the above deformed models. Moreover, \cite{BGRao09,BGRao10,CDFF10} show us that free probability can also allow to locate 
isolated  eigenvalues of deformed matricial models.
In particular, in  \cite{CDFF10}, the authors point out that the subordination function relative to the free additive convolution provides a good understanding of the outliers of deformed Wigner matrices.
 We will see in this paper that  the subordination function relative to the free (additive or multiplicative) convolution
 plays again an important part in the asymptotic behaviour of the  eigenvectors relative to the outliers.  We introduce in the following section some  results concerning free convolution that will be fundamental later on.

\section{ Free convolution}\label{freeconv}
Free convolutions appear as natural analogues of the classical convolutions in the context of free probability theory.
Denote by ${\cal M}$ the set of probability measures supported on the real line and by ${\cal M}^+$ the ones supported on $[0; + \infty[$. For $\mu$ and $\nu$ in ${\cal M}$ one defines the free additive convolution $\mu \boxplus \nu$ of $\mu$ and $\nu$ as the distribution of $X+Y$ where $X$ and $Y$ are free self adjoint random variables with distribution $\mu$ and $\nu$. For $\mu$ and $\nu$ in ${\cal M}^+$,  the free multiplicative convolution $\mu \boxtimes \nu$ of $\mu$ and $\nu$ is the distribution of $X^{\frac{1}{2}}YX^{\frac{1}{2}}$ where $X$ and $Y$ are free positive random variables with distribution $\mu$ and $\nu$. We refer the reader to \cite{VDN92} for an introduction to free probability theory
and to \cite{Voiculescu86,Voiculescu87} and \cite{BercoviciVoiculescu} for free convolutions.
In this section, we recall the  analytic approach developped in \cite{Voiculescu86,Voiculescu87} to calculate the free convolutions of measures, we present
the important subordination property, and describe more deeply subordination functions relative to free additive  convolution by a semi-circular distribution and free multiplicative convolution by a Marchenko-Pastur distribution.
We also recall  characterizations of the complement of the support of these convolutions.
\subsection{Additive Free convolution}
Let $\tau $ be a probability measure on $\R$. 
Its Stieltjes transform $g_\tau $ is analytic on the complex upper half-plane $\C^+$. 
There exists a domain $$D_{\alpha , \beta } = \{ u+iv \in \C, |u| < \alpha v, v > \beta \}$$
on which $g_\tau $ is univalent. 
Let $K_\tau $ be its inverse function, defined on $g_\tau (D_{\alpha , \beta })$, and 
$$R_\tau (z) = K_\tau (z) - \frac{1}{z}.$$
Given two probability measures $\tau $ and $\nu $, 
there exists a unique probability measure $\lambda $ such that
$$R_\lambda = R_\tau + R_\nu $$
on a domain where these functions are defined. 
The probability measure $\lambda $ is called 
the additive free convolution of $\tau $ and $\nu $ and denoted by $\tau \boxplus \nu $. 
\subsubsection{Subordination property}
The  free additive  convolution of probability measures has an important property, 
called subordination, which can be stated as follows.
\begin{proposition} 
let $\tau $ and $\nu $ be two probability measures on $\R$; 
there exists a unique analytic map $F^{(a)}: \C^+ \rightarrow \C^+$
such that \begin{equation} \label{subeq}\forall z \in \C^+ ,  ~~~~g_{\tau \boxplus \nu}(z)= g_\nu (F^{(a)}(z)),\end{equation}
$$\, F^{(a)}(\overline{z})=\overline{F^{(a)}(z)}, \, \Im F^{(a)}(z)\geq \Im z, \, \lim_{y \rightarrow + \infty} \frac{F^{(a)}(iy)}{iy}=1.$$
\end{proposition}
This phenomenon was first observed by D. Voiculescu under a genericity assumption in \cite{Voiculescu93}, 
and then proved in generality in \cite{Biane98} Theorem 3.1. 
Later, a new proof of this result was given in \cite{BelBer07}, 
using a fixed point theorem for analytic self-maps of the upper half-plane. 

\subsubsection{Free convolution by a semicircular distribution}

In \cite{Biane97b}, P. Biane provides a deep study of the free additive convolution by a semicircular distribution. 
We first recall here some of his results that will be useful in our approach. Let $\nu $ be a probability measure on $\R$. 
When $\tau$ in (\ref{subeq}) is the  semi-circular distribution $\mu_\sigma$, let us denote by $F^{(a)}_{\sigma,\nu}$ the subordination function. In \cite{Biane97b},
P. Biane  introduces the set 
\begin{equation}\label{Ome}\Omega _{\sigma , \nu }:=\{ u+iv \in \C^+, v > v_{\sigma , \nu }(u)\},\end{equation}
where the function $v_{\sigma , \nu }: \R \rightarrow \R^+$ is defined by 
$$v_{\sigma , \nu }(u) = \inf \left\{v \geq 0, \int_{\R} \frac{d\nu (x)}{(u-x)^2+v^2} \leq \frac{1}{\sigma ^2}\right\}.$$ 
The boundary of $\Omega _{\sigma , \nu }$ is the graph of the continuous function $v_{\sigma , \nu }$.
P. Biane proves the following

\begin{proposition}\cite{Biane97b}\label{homeo} 
The map 
\begin{equation}\label{defdeH} H_{\sigma , \nu }: z \longmapsto z+\sigma ^2 g_\nu (z)\end{equation}
is a homeomorphism from $\overline{\Omega _{\sigma , \nu }}$ to $\C^+ \cup \R$ 
which is conformal from $\Omega _{\sigma , \nu }$ onto $\C^+$. 
$F^{(a)}_{\sigma , \nu }: \left\{\begin{array}{ll} \C^+ \cup \R \rightarrow \overline{\Omega _{\sigma ,\nu }}\\ 
z \rightarrow z-\sigma ^2g_{\mu _\sigma \boxplus \nu }(z) \end{array}\right.$ 
is the inverse function of $H_{\sigma , \nu }$.
\end{proposition}

\noindent Considering $ H_{\sigma , \nu }$ as an analytic map defined in the whole upper half-plane $\C^+$, 
it can be easily seen  that \begin{equation}\label{imageinverse}\Omega _{\sigma , \nu }=( H_{\sigma , \nu })^{-1}(\C^+). \end{equation}
\noindent {\bf In the following, we will denote $H_{\sigma , \nu }$ by $H$ to simplify the writing.}
\begin{remarque}\label{imagedeF}
Note that according to Proposition \ref{homeo}, 
$$F^{(a)}_{\sigma , \nu }(\mathbb{R}) =\partial \Omega_{\sigma , \nu } =\{u+i v_{\sigma , \nu }(u), u\in \mathbb{R}\}$$
\noindent so that we have the following equivalence
$$u\in  F^{(a)}_{\sigma , \nu }(\mathbb{R}) \cap \mathbb{R} \Longleftrightarrow v_{\sigma , \nu }(u)=0.$$
\noindent The following characterization of the elements of the complement of the support of $\nu$ which are in the image of $\mathbb{R}$ by $F^{(a)}_{\sigma , \nu }$ readily follows:
\begin{equation}\label{imagereelle}u \in  F^{(a)}_{\sigma , \nu }(\mathbb{R})  \cap      \mathbb{R}\setminus {\rm supp~} (\nu) \Longleftrightarrow
u \in ^c{\rm supp~} (\nu) , H'(u) \geq 0.\end{equation}
\end{remarque}
\noindent 
In \cite{Biane97b}, P. Biane  obtains a description of the support of $\mu _\sigma \boxplus \nu $ from which,
when $\nu $ is a compactly supported probability measure, 
the authors deduce in \cite{CDFF10} a characterization of the complement of the support 
of $\mu _\sigma \boxplus \nu $ involving the support of $\nu $ and $H$.

\begin{proposition}\label{Caract}
$$x \in  ^c {\rm supp}(\mu _\sigma \boxplus \nu ) \Leftrightarrow 
\exists u \in {\cal O}^{(a)} {\rm ~such~that~} x=H(u)$$
where ${\cal O}^{(a)}$ is the open set
\begin{eqnarray} {\cal O}^{(a)}&:=&\left\{ u \in ^c {\rm supp}(\nu ),  \,  \, H'(u) > 0\right\} \nonumber \\& =&\left\{ u \in ^c {\rm supp}( \nu),\,  \sigma^2 \int \frac{1}{(u-t)^2}d \nu(t) <1\right\} \label{defThetaW}.\end{eqnarray}
\end{proposition}  
\begin{remarque}\label{oainclu}  ${\cal O}^{(a)} \subset \partial \Omega_{\sigma , \nu }$
\end{remarque}
\noindent This readily follows from Remark \ref{imagedeF}.
\begin{remarque}\label{Hcroit}
Note that 
if $u_1 < u_2$ are in $\left\{ u \in ^c {\rm supp~} (\nu), H'(u)\geq 0 \right\}$, one has 
$H(u_1)\leq H(u_2).$
 Indeed, by Cauchy-Schwarz inequality, we have 
\begin{eqnarray*}
H(u_2)-H(u_1) & = & (u_2-u_1) \bigg [ 1 - \sigma ^2 \int_{\R}\frac{d\nu (x)}{(u_1-x)(u_2-x)} \bigg ] \\
& \geq & (u_2-u_1) \bigg [ 1 - \sigma ^2 \sqrt{ (-g'_{\nu }(u_1))(-g'_{\nu }(u_2))} \bigg ] \geq 0.
\end{eqnarray*}
\end{remarque}

\subsection{Multiplicative free convolution}

Let $\tau \neq \delta_0$ be a probability measure on $[0;+ \infty[$. Define the analytic function $$\Psi_\tau (z)=\int \frac{tz}{1-tz}d\tau(t)=\frac{1}{z}g_\tau(\frac{1}{z}) -1,$$
\noindent for complex values of $z$ such that $\frac{1}{z}$ is not in the support of $\tau$.
$\Psi_\tau$ determines uniquely the measure $\tau$ and it is univalent in the left half-plane $\{z \in \mathbb{C},\, \Re z <0\}$.

\noindent Then one may determine an analytic function $S_\tau$ such that 
$$\Psi_\tau \left[ \frac{z}{z+1} S_\tau(z) \right]=z$$
\noindent in some domain (which will contain at least some interval to the left of zero) and then 
$S_{\mu \boxtimes \nu}=S_\mu S_\nu.$ (see \cite{Voiculescu87}).

\subsubsection{Subordination property}
Free multiplicative convolution also presents a subordination phenomenon first proved in \cite{Biane98} (see also \cite{BelBer07}).
\begin{proposition}\label{submult}
Let $\tau \neq \delta_0$ and $\nu \neq \delta_0$ be two probability measures on $[0;+\infty[$. There exists a unique analytic
map $F^{(m)}_{\tau,\nu}$  defined on $\mathbb{C}\setminus [0;+\infty[$ such that
\begin{equation}\label{eqsubmult}\forall \, z \in \mathbb{C}\setminus [0;+\infty[, \, \Psi_{\nu\boxtimes \tau}(z) =\Psi_\nu(F^{(m)}_{\tau,\nu}(z))\end{equation}
and $$\forall \, z \in \mathbb{C}^+, \; F^{(m)}_{\tau,\nu}(z) \in \mathbb{C}^+,\, F^{(m)}_{\tau,\nu}(\overline{z})
=\overline{F^{(m)}_{\tau,\nu}(z)},\, 
\, \arg (F^{(m)}_{\tau,\nu}(z))\geq \arg(z).$$
\end{proposition}
\subsubsection{Multiplicative free convolution with a Marchenko-Pastur distribution}
Let us determinate the subordination function relative to the free multiplication by  a Marchenko-Pastur distribution.
We can deduce from (\ref{MPeq}) that for any $z \in \mathbb{C}^+$,
\begin{equation}\label{MPeq2} g_{ \mu_{\mbox{\tiny{MP}},c}\boxtimes \nu}(z)=\int \frac{1}{z- t(1-c +cz g_{ \mu_{\mbox{\tiny{MP}},{c}}\boxtimes \nu}(z))}d\nu(t)\end{equation}
and then that $\forall z \in \mathbb{C}\setminus [0;+\infty[$,
\begin{equation}\label{psi}\Psi_{ \mu_{\mbox{\tiny{MP}},c}\boxtimes \nu} (z)=\Psi_\nu \left(z-cz +cg_{ \mu_{\mbox{\tiny{MP}},c}\boxtimes \nu}(\frac{1}{z})\right).\end{equation}
\noindent Note that $$z-cz +cg_{ \mu_{\mbox{\tiny{MP}},c}\boxtimes \nu}(\frac{1}{z})=g_{\tau_{c,\nu}}(\frac{1}{z})$$
where {\bf $\tau_{c,\nu}$ is the limiting spectral distribution of $\frac{1}{p} B_N^* A_N B_N$}.\\

\noindent It is clear that $\forall z \in \mathbb{C}^+, \; g_{\tau_{c,\nu}}(\frac{1}{z}) \in \mathbb{C}^+,\,  g_{\tau_{c,\nu}}(\frac{1}{\overline{z}})
=\overline{ g_{\tau_{c,\nu}}(\frac{1}{z})}.$ Moreover, since, using (\ref{MPeq2}), we have 
$$g_{\tau_{c,\nu}}(\frac{1}{z})=z \left[(1-c) +c
\int \frac{1}{1- tg_{\tau_{c,\nu}}(\frac{1}{z})}d\nu(t)\right],$$
\noindent it is easy to see that
$\arg (g_{\tau_{c,\nu}}(\frac{1}{z}))\geq \arg(z).$\\
Therefore,
 denoting by $F^{(m)}_{c,\nu}$ the subordination function in (\ref{eqsubmult}) when $\tau= \mu_{\mbox{\tiny{MP}},c}$, we have that 
 
$$ F^{(m)}_{c,\nu}(z)= z-cz +cg_{ \mu_{\mbox{\tiny{MP}},c}\boxtimes \nu}(\frac{1}{z})=g_{\tau_{c,\nu}}(\frac{1}{z}).$$
Now, we are going to present the characterization of the complement of the support of 
$\tau_{c,\nu}$ provided by  Choi and Silverstein  in \cite{CHOISILVER} .  Note that the supports of $ \mu_{\mbox{\tiny{MP}},c}\boxtimes \nu$ and  $\tau_{c,\nu}$ obviously coincide on $]0;+\infty[$.\\
 ~~

\noindent According to \cite{BaiSil06}  p 113, for each $z \in \mathbb{C^+}$, $g_{\tau_{c,\nu}}(z)$ is the unique solution $Z$ in $\mathbb{C}^{-}$
of the equation \begin{equation} Z =\frac{1}{z-c\int \frac{t}{1-tZ}d\nu(t)}\end{equation}
so that \begin{equation}\label{composition} z= {\cal Z}_{c,\nu}(g_{\tau_{c,\nu}}(z))\end{equation}
where

\begin{equation}\label{defzcnu}{\cal Z}_{c,\nu}(x)=\frac{1}{x}+c \int \frac{t}{1-tx}d\nu(t).\end{equation}

\noindent {\bf In the following, we will denote ${\cal Z}_{c,\nu}$ by ${\cal Z}$ to simplify the writing.}
\begin{proposition}\label{ChoiSilver}\cite{CHOISILVER}
If $u \in ^c {\rm supp}( \tau_{c,\nu})$, then $s=g_{\tau_{c,\nu}}(u)$ satisfies 
\begin{enumerate}
\item $s \in \mathbb{R}\backslash\{0\},$
\item $\frac{1}{s} \in ^c {\rm supp}( \nu)$,
\item ${\cal Z}{'}(s)<0$.
\end{enumerate}
Conversely if $s$ satisfies (1), (2) and (3), then $u = {\cal Z}(s)\in ^c {\rm supp}( \tau_{c,\nu})$.

\end{proposition}

\noindent In particular, letting $z$ converge towards any element $u$ of $^c \mbox{supp}( \tau_{c,\nu})$ in (\ref{composition})
leads to \begin{equation}\label{compodansunsens}u= {\cal Z}(g_{\tau_{c,\nu}}(u)).\end{equation}

\noindent In \cite{CHOISILVER}, the authors proved also that $\lim_{\begin{array}{ll}z\rightarrow u\\z \in \mathbb{C}^+ \end{array}} g_{\tau_{c,\nu}}(z)=g_{\tau_{c,\nu}}(u)$ exists for any $u$ in $\mathbb{R}\setminus\{0\}$. 

\noindent We include here for the convenience of the reader some basic facts that will be used later on.
\begin{remarque}\label{compositionlimite}
For any $u \in \R\setminus \{0\}$ such that $g_{\tau_{c,\nu}}(u)\in \mathbb{R}\setminus\{0\}$ and $\frac{1}{g_{\tau_{c,\nu}}(u)} \in ^c {\rm supp~}(\nu)$, we have ${\cal Z}[g_{\tau_{c,\nu}}(u)]=u$. 
\end{remarque}
\noindent This readily follows by letting $z$ goes to $u$ in (\ref{composition}).\\

\noindent Let us introduce the open set
\begin{eqnarray} \label{defThetaS} {\cal O}^{(m)}&:=&\left\{ u \in ^c \mbox{supp}( \nu)\setminus \{0\},\, {\cal Z}^{'}(\frac{1}{u})<0\right\} \nonumber \\& =&\left\{ u \in  ^c \mbox{supp}( \nu)\setminus \{0\},\,  c \int \frac{t^2}{(u-t)^2}d \nu(t) <1\right\} \label{defThetaS}.\end{eqnarray}

\begin{remarque}\label{CScomp}
For any $u$ in ${\cal O}^{(m)}$, $g_{\tau_{c,\nu}}[{\cal Z}(\frac{1}{u})]=\frac{1}{u}$.
\end{remarque}
\noindent Let us prove Remark \ref{CScomp}. According to Proposition \ref{ChoiSilver}, for any $u \in {\cal O}^{(m)}$,
${\cal Z}(\frac{1}{u}) \in ^c \mbox{supp}( \tau_{c,\nu})$ and then according to the same Proposition \ref{ChoiSilver},
$\frac{1}{g_{\tau_{c,\nu}}[{\cal Z}(\frac{1}{u})]}$ also belongs to ${\cal O}^{(m)}$. Now, for any $a\neq b$ in ${\cal O}^{(m)}$,
\begin{equation}\label{croi} {\cal Z}(\frac{1}{b})-{\cal Z}(\frac{1}{a})={(b-a)} \left[ 1-c \int \frac{t^2}{({a}-t)({b}-t)} d \nu(t) \right];\end{equation}
\noindent by Cauchy-Schwartz inequality, 
\begin{eqnarray}\hspace*{-5mm}\left| c \int \frac{t^2}{({a}-t)({b}-t)} d \nu(t) \right| &\leq & \left\{c\int \frac{t^2}{({a}-t)^2} d \nu(t)\right\}^{\frac{1}{2}}  \left\{c\int \frac{t^2}{(b-t)^2} d \nu(t)\right\}^{\frac{1}{2}} \label{borne}\\\hspace*{-5mm}& < &1.
\nonumber
\end{eqnarray}
 Hence we can conclude that
\begin{equation} \label{difference}  \mbox{~ for any~} a\neq b\;  \mbox{~in~} {\cal O}^{(m)}, \, {\cal Z}(\frac{1}{b}) \neq {\cal Z}(\frac{1}{a}).\end{equation}
Since using (\ref{compodansunsens}) we have ${\cal Z}[g_{\tau_{c,\nu}}[{\cal Z}(\frac{1}{u})]]={\cal Z}(\frac{1}{u})$, we can then deduce that $g_{\tau_{c,\nu}}[{\cal Z}(\frac{1}{u})]=\frac{1}{u}.$ $\Box$\\
~~

\begin{remarque}\label{croissance globaleS} Using  (\ref{croi}) and (\ref{borne}), we have 
for any $a \leq b$ in the set $\left\{ u \neq 0,\, u \in  ^c {\rm supp}( \nu),\,  {\cal Z}^{'}(\frac{1}{u})\leq 0\right\}$ that  ${\cal Z}(\frac{1}{a}) \leq {\cal Z}(\frac{1}{b})$.
\end{remarque}

\begin{remarque}\label{nonnul}
For any $u \in {\cal O}^{(m)}\cap ]0;+\infty[$, we have $ {\cal Z}(\frac{1}{u})>0$. 
\end{remarque}
\noindent Indeed, assume that $ {\cal Z}(\frac{1}{u}) \leq 0$.
According to Proposition \ref{ChoiSilver}, ${\cal Z}(\frac{1}{u})\in ^c \mbox{supp}( \tau_{c,\nu})$;   $ {\cal Z}(\frac{1}{u}) \leq~0$ implies that  ${\cal Z}(\frac{1}{u})$ is on the left hand side of $ \mbox{supp}( \tau_{c,\nu})$ and therefore that 
$g_{\tau_{c,\nu}}({\cal Z}(\frac{1}{u}))\leq 0$. This leads to a contradiction with  Remark \ref{CScomp} saying that $g_{\tau_{c,\nu}}({\cal Z}(\frac{1}{u}))=\frac{1}{u}>0.$
\begin{remarque}\label{pasdansimagedeg}
For any $u\neq 0$ in $^c {\rm supp~}( \nu)$ such that ${\cal Z}^{'}(\frac{1}{u})> 0$, we have $\frac{1}{u} \notin 
g_{\tau_{c,\nu}} (\mathbb{R}\setminus\{0\})$.
\end{remarque}
\noindent Indeed, let us assume that there exists $v \in \mathbb{R}\setminus\{0\}$ such that $\frac{1}{u} =
g_{\tau_{c,\nu}}(v)$. According to Proposition \ref{ChoiSilver}, $v$ belongs to $\mbox{supp}( \tau_{c,\nu})$.
Using Remark \ref{compositionlimite}, since  $g_{\tau_{c,\nu}}(v)\in \mathbb{R}\setminus\{0\}$ and $\frac{1}{g_{\tau_{c,\nu}}(v)} \in ^c \mbox{supp~}(\nu)$, we have ${\cal Z}[g_{\tau_{c,\nu}}(v)]=v$. It follows that for any $y>0$,
\begin{eqnarray}
1&=& \frac{{\cal Z}[g_{\tau_{c,\nu}}(v+iy)]-{\cal Z}[g_{\tau_{c,\nu}}(v)]}{iy} \nonumber\\
&=&
\frac{{\cal Z}[g_{\tau_{c,\nu}}(v+iy)]-{\cal Z}[g_{\tau_{c,\nu}}(v)]}{g_{\tau_{c,\nu}}(v+iy)-g_{\tau_{c,\nu}}(v)}      \times  \frac{g_{\tau_{c,\nu}}(v+iy)-g_{\tau_{c,\nu}}(v)} {iy} .\label{fractions}
\end{eqnarray}
Since $g_{\tau_{c,\nu}}(v+iy)$ converges towards $g_{\tau_{c,\nu}}(v) =\frac{1}{u}$ and $ {\cal Z} $ is holomorphic in a neighborhood of $g_{\tau_{c,\nu}}(v)$, letting $y $ tends to zero the first factor on the right hand side of (\ref{fractions}) converges towards ${\cal Z}^{'} (g_{\tau_{c,\nu}}(v))= {\cal Z}^{'}(\frac{1}{u})>0.$ This implies that $\Re\left[ \frac{g_{\tau_{c,\nu}}(v+iy)-g_{\tau_{c,\nu}}(v)} {iy}\right]$ converges towards $\frac{1}{{\cal Z}^{'}(\frac{1}{u})}>0$ when $y$ tends to zero. Now, for any $y>0$, we have 
\begin{eqnarray*}\Re\left[ \frac{g_{\tau_{c,\nu}}(v+iy)-g_{\tau_{c,\nu}}(v)} {iy}\right]&=& 
\frac{\Im \left[ g_{\tau_{c,\nu}}(v+iy)-g_{\tau_{c,\nu}}(v)\right]} {y}\\
&=&\frac{\Im g_{\tau_{c,\nu}}(v+iy)}{y} <0
\end{eqnarray*}
\noindent which leads to a contradiction. $\Box$

\section{Main results}

As noticed in the previous section, we have the  following characterization of the complement of the support of the limiting spectral distribution of $M_N^W$.

\begin{equation}\label{car1}x \in  ^c{\rm supp}(\mu _\sigma \boxplus \nu ) \Leftrightarrow 
\exists u \in {\cal O}^{(a)} {\rm ~such~that~} x=H(u) .\end{equation}
Moreover, we can deduce, using Proposition \ref{homeo} and Remark \ref{oainclu}, that :\\ 
 $x \mapsto F^{(a)}_{\sigma,\nu}\left(x\right)$ is a bijection from $ ^c{\rm supp}(\mu _\sigma \boxplus \nu )$ onto ${\cal O}^{(a)}$ with inverse $H$.\\

\noindent Since the supports of $ \mu_{\mbox{\tiny{MP}},c}\boxtimes \nu$ and  $\tau_{c,\nu}$  coincide on $]0;+ \infty[$, we can also deduce from the previous section  the  following characterization of the restriction to $\mathbb{R}\setminus\{0\}$ of the complement of the  support of the limiting spectral distribution of $M_N^S$:
\begin{equation}\label{car2}x \neq 0, x \in  ^c \mbox{supp}(  \mu_{\mbox{\tiny{MP}},c}\boxtimes \nu) \Leftrightarrow \exists u \in {\cal O}^{(m)}, {\cal Z} (\frac{1}{u}) \neq 0, {\rm ~such~that~}
x={\cal Z} (\frac{1}{u}).\end{equation}
\noindent Moreover, we can deduce from Remark \ref{compositionlimite} and Remark \ref{CScomp} that:\\ 
 $x \mapsto F^{(m)}_{c,\nu}\left(\frac{1}{x}\right) (=g_{\tau_{c,\nu}}(x))$ is a bijection from $ ^c \mbox{supp}(  \mu_{\mbox{\tiny{MP}},c}\boxtimes \nu)\setminus \{0\}$ onto the set $\left\{\frac{1}{u}, u \in {\cal O}^{(m)}, {\cal Z} (\frac{1}{u})\neq 0\right\}$ with inverse ${\cal Z}$.\\
Note that the limiting  mass at zero was studied in \cite{CHOISILVER}:
$$ \mu_{\mbox{\tiny{MP}},c}\boxtimes \nu (0) =\left\{\begin{array}{ll}\nu(0) \mbox{~if~} c(1-\nu(0))\leq 1,\\
1- \frac{1}{c} \mbox{~if~} c(1-\nu(0))> 1.\end{array} \right.$$

Actually, according to \cite{CDFF10} (resp. \cite{RaoSil09,BaiYao08b}),
 the spikes $\theta_j$ in $\Theta=\{\theta_1;\ldots;\theta_J\}$ of the perturbation matrix $A_N$  that will generate eigenvalues of $M_N^W$ (resp. $M_N^S$) which deviate from the bulk are exactly those belonging to ${\cal O}^{(a)}$ (resp. ${\cal O}^{(m)}$)
 and the corresponding limiting points 
outside the support of $\mu _\sigma \boxplus \nu $ (resp. $\tau_{c,\nu}$) will be given by
$H(\theta_j)$ (resp. ${\cal Z} (\frac{1}{\theta_j})$). 
Note that the results in \cite{RaoSil09,BaiYao08b} are not formulated in that way since the authors  do not 
deal with subordination function but as already mentioned we choose to express all the results 
using these functions $H$ and ${\cal Z}$ related to the subordination functions  for further generalizations.
Hence adopting the notations of the first column of the following array standing for both the corresponding elements of the second column (deformed Wigner matrix case) and the third column (sample covariance matrix case),
we can present 
a common  formulation of these results. 
\begin{equation}\label{correspondance}\mbox{
\begin{tabular}{|c|c|c|} \hline
$M_N$& $M_N^W$ &  $ M_N^S$\\ \hline
$\mu_{\mbox{\tiny{LSD}}}$& $ \mu _\sigma \boxplus \nu$ & $  \mu_{\mbox{\tiny{MP}},c}\boxtimes \nu$\\ \hline
${\cal O} $&${\cal O}^{(a)}$ &${\cal O}^{(m)}$\\ \hline
$ \Theta_o$ & $\Theta \cap {\cal O}^{(a)}=\{\theta_i \in \Theta, H' (\theta_j)>0\}$ &$ \Theta \cap {\cal O}^{(m)}=\{\theta_i \in \Theta, {\cal Z}^{'}(\frac{1}{\theta_j})<0\}  $\\ \hline
$ \rho _{\theta _j}$& $H(\theta_j)$& ${\cal Z}(\frac{1}{\theta_j}) $\\ \hline
\end{tabular}}
 \end{equation}
\begin{theoreme}{\label{ThmASCV}}\cite{RaoSil09,BaiYao08b,CDFF10}
Let  $\theta_j$ be  in $ \Theta_o$ 
and  denote by $n_{j-1}+1, \ldots , n_{j-1}+k_j$ the descending ranks of $\theta _j$ among the eigenvalues of $A_N$.
Then the $k_j$ eigenvalues $(\lambda_{n_{j-1}+i}(M_N), \, 1 \leq i \leq k_j)$ 
converge almost surely outside the support of $\mu_{\mbox{\tiny{LSD}}}$
towards $\rho _{\theta _j}$.
Moreover, these eigenvalues asymptotically separate from the rest of the spectrum since 
(with the conventions that $\lambda_0(M_N)=+\infty$ and $\lambda_{N+1}(M_N)=-\infty$)
there exists  $0< \delta_0 $ such that 
\noindent almost surely for all large N, \begin{equation}\label{sepraj}\lambda_{n_{j-1}}(M_N) > \rho _{\theta _j} + \delta_0 \, \mbox{~and~} \, \lambda_{n_{j-1}+k_j +1}(M_N) < \rho _{\theta _j} - \delta_0 
.\end{equation}
\end{theoreme}

 The aim of this paper is  to study  how the corresponding eigenvectors of the deformed model project onto those of the perturbation.
Here is  the main result of the paper still adopting the notations of the first column of the  array (\ref{correspondance}) in order to present   
a unified approach.
\begin{theoreme}\label{cvev1p1}
 Let  $\theta_j$ be in $\Theta_o$
and denote by $n_{j-1}+1, \ldots , n_{j-1}+k_j$ the descending ranks of $\theta _j$ among the eigenvalues of $A_N$.
 Let    
 $\xi(j)$ be  a normalized  eigenvector of $M_N$  relative to one of the eigenvalues  $(\lambda_{n_{j-1}+q}(M_N)$, $1\leq q\leq k_j)$.
Then,
when $N$ goes to infinity,
\begin{itemize}
\item[(i)] $\left\| P_{\mbox{Ker~}(\theta_j I_N-A_N)}\xi(j)\right\|^2_2  \stackrel{a.s}{\rightarrow} \tau(\theta_j)$\\
 ~~

\noindent where  \begin{equation}\label{tau}\tau(\theta_j)=\left\{\begin{array}{ll} H'(\theta_j)  \mbox{~if~} M_N=M_N^W,\\
- \frac{{\cal Z}^{'}(\frac{1}{\theta_j})}{\theta_j {\cal Z}(\frac{1}{\theta_j}) } \mbox{~if~} M_N=M_N^S.
 \end{array}\right.\end{equation}
\noindent  Note that we have explicitly $$H'(\theta_j) = 1 -\sigma^2 \int \frac{1}{(\theta_j-x)^2}d\nu(x),$$ 
 $$- \frac{{\cal Z}'(\frac{1}{\theta_j})}{\theta_j {\cal Z}(\frac{1}{\theta_j}) }=
 \frac{1- c\int \frac{x^2}{(\theta_j-x)^2}d\nu(x)}{1+c \int \frac{x}{(\theta_j-x)}d\nu(x)}.$$
\item[(ii)] for any  $\theta_l$ in $\Theta\setminus\{\theta_j\}$, 
$$\left\| P_{\mbox{Ker~}(\theta_l I_N-A_N)}\xi(j)\right\|_2  \stackrel{a.s}{\rightarrow} 0.$$

\end{itemize}
\end{theoreme}

\noindent {\bf Example:} Let us consider the perturbation matrix
$$A_N  =  \displaystyle{{\rm{ diag
}}(2, \frac{3}{2},0,\underbrace{-1,\ldots
  ,-1}_{\frac{N}{2}},   \underbrace{1,\ldots,
 1 }_{\frac{N}{2}-3})},$$ whose limiting spectral distribution is 
$\nu =\frac{1}{2}\delta_1 + \frac{1}{2}\delta_{-1}$. Thus, the set of the spikes of $A_N$ is $\Theta=\{2; \frac{3}{2};0\}$.
Let us consider the corresponding  deformed Wigner model assuming moreover that $\sigma^2=\frac{1}{2}$. Then,
 $H(u)=u+\frac{1}{4}\frac{1}{(u-1)}+\frac{1}{4}\frac{1}{(u+1)}$. One can check that  the support of $\mu_{\frac{1}{\sqrt{2}}}\boxplus \nu$ has two connected components which are symmetric with respect to zero.
Since  $H'(2)=\frac{13}{18}>0$ and $2$ is the largest eigenvalue of $A_N$, according to Theorem \ref{ThmASCV}, when $N$ goes to infinity,
the largest eigenvalue of the deformed Wigner model $M_N^W$  converges almost surely towards $H(2)=\frac{7}{3}$ (on the right hand side of the support of 
$\mu_{\frac{1}{\sqrt{2}}}\boxplus \nu$). 
 Note that, since $H'(\frac{3}{2})<0$, the second largest  eigenvalue of $M_N^W$ sticks to the bulk.
Moreover, since $H'(0)=\frac{1}{2}>0$ and the descending rank of $0$ among the eigenvalues of $A_N$ is $\frac{N}{2}$,
according to Theorem \ref{ThmASCV}, when $N$ goes to infinity, $\lambda_{\frac{N}{2}}(M_N^W)$ converges almost surely towards $H(0)=0$ which is
between the two connected components of the support of $\mu_{\frac{1}{\sqrt{2}}}\boxplus \nu$.
Now, denote by $\{e_1,\ldots,e_N\}$ the canonical basis of $\mathbb{C}^N$.
Since $e_1$ is an eigenvector relative to $2$, $e_2$ is an eigenvector relative to $\frac{3}{2}$ and $e_{3}$ is an eigenvector relative to $0$, according to Theorem \ref{cvev1p1}, if $\xi$ denotes a normalized eigenvector associated to the largest eigenvalue of $M_N^W$, $^t \xi=\left( \xi^{(1)},\ldots,\xi^{(N)}\right)$, then 
$\vert\xi^{(1)}\vert  \stackrel{a.s}{\rightarrow} \sqrt{H'(2)}=\frac{\sqrt{13}}{3\sqrt{2}}$ and, for $i=2,3$, $\vert\xi^{(i)}\vert  \stackrel{a.s}{\rightarrow} 0$
when $N$ goes to infinity. Similarly,  if $\xi$ denotes a normalized eigenvector associated to  $\lambda_{\frac{N}{2}}(M_N^W)$, $^t \xi=\left( \xi^{(1)},\ldots,\xi^{(N)}\right)$, then 
$\vert\xi^{(3)}\vert  \stackrel{a.s}{\rightarrow} \sqrt{H'(0)}=\frac{1}{\sqrt{2}}$ and, for $i=1,2$, $\vert\xi^{(i)}\vert  \stackrel{a.s}{\rightarrow} 0$
when $N$ goes to infinity. \\


Actually, in order to establish Theorem \ref{cvev1p1}, we will first prove Proposition \ref{cvev} below since
when $k_j \neq 1$, the method used in this paper does not allow us  to tackle directly the orthogonal projection of each eigenvector separately to prove (i).
\begin{proposition}\label{cvev} Let  $\theta_j$ be in $\Theta_o$
and denote by $n_{j-1}+1, \ldots , n_{j-1}+k_j$ the descending ranks of $\theta _j$ among the eigenvalues of $A_N$.
Denote by  
 $\xi_1(j),\ldots, \xi_{k_j}(j)$ an orthonormal system of eigenvectors associated to $(\lambda_{n_{j-1}+n}(M_N)$, $1\leq n\leq k_j)$.
Then, for any  $\theta_l $ in $\Theta$, 
 when $N$ goes to infinity, $$\sum_{n=1}^{k_j}\left\| P_{\mbox{Ker~}(\theta_l I_N-A_N)}\xi_n(j)\right\|^2_2  \stackrel{a.s}{\rightarrow} \delta_{l,j}k_j\tau(\theta_j)$$
where $\tau(\theta_j)$ is defined by (\ref{tau}).
\end{proposition}
In Section \ref{redu}, we will explain how we can 
  deduce (i) of Theorem \ref{cvev1p1} from Proposition \ref{cvev} using a perturbation trick and ideas of \cite{EI}.\\
~~

\section{Reduction of the proof of (i) of Theorem \ref{cvev1p1}  to the case  of a spike with multiplicity one}\label{redu}

Note first that, dealing with a spike $\theta_j$ in $\Theta_o$ with multiplicity one, the statements of   Theorem \ref{cvev1p1} and Proposition \ref{cvev}
are equivalent. Thus, in this section, we show how to deduce (i) of Theorem \ref{cvev1p1} dealing with a spike $\theta_j$ with multiplicity $k_j \neq 1$
from the hypothesis that (i) is true dealing with a spike with multiplicity one.
We will need the following lemmas.
\begin{lemme}\label{Proj}
Let ${\cal V}$ be a vector subspace of $\mathbb{C}^N$ with an orthonormal basis $V_1,\ldots,V_k$.
Let $\alpha$ be  in $[0;+\infty[$.
For any $m=1,\ldots,k$, let $\alpha_m$  be in $[0;+\infty[$.
 Let $E$ be a vector subspace of $\mathbb{C}^N$ with an orthonormal basis $\xi_1,\ldots,\xi_k$.
Then, there exists a sequence $a_N \geq 0$, depending on the $ \vert \langle V_i,\xi_n\rangle\vert$, $ \alpha_m$, $i,n,m \in \{1,\ldots,k\}^3$,
 such that, for any vector $u$ in the unit sphere of $ \mathbb{C}^N$,
$$\left|  \left\|P_{\cal V} u\right\|_2^2-\alpha\right| \leq (2k+\alpha) \left\|P_{{E}^{\bot}} u \right\|_2 +\max_{m=1}^{k} \left| \alpha_m-\alpha \right| +a_N,$$
\noindent and, if,  for any $ m,n$ in $\{1,\ldots,k\}^2$, \begin{equation}\label{hypconv} \vert \langle V_m,\xi_n\rangle\vert^2 
\rightarrow \delta_{m,n} \alpha_m\mbox{~~~when $N$ goes to infinity},
\end{equation}
then $a_N $ converges to zero when $N$ goes to infinity.
\end{lemme}
\noindent {\bf Proof:} Throughout the proof, we will often use the following  obvious inequalities without mentioning them:\\
 for any vectors $u_1$ and $u_2$ in the unit sphere of $ \mathbb{C}^N$,
 $\left| \langle u_1u_2\rangle \right| \leq  \left\|  u_1 \right\|_2 \left\| u_2 \right\|_2=1$,
$\left| \langle P_{E}u_1,u_2\rangle \right|\leq  \left\| P_{E} u_1 \right\|_2  \left\|u_2\right\|_2  \leq  \left\|  u_1 \right\|_2 \left\| u_2 \right\|_2=1.$\\
 ~~

\noindent  We have for each $m$ in $\{1,\ldots,k\}$, for any vector $u$ in the unit sphere of $ \mathbb{C}^N$,\\

~~~~$\left| ~\left| \langle u,V_m\rangle \right|^2 -\left| \langle P_{E} u,V_m\rangle \right|^2 ~ \right|$
\begin{eqnarray} &=&\left| \left( \left| \langle u,V_m\rangle \right| -\left| \langle P_{E}u,V_m\rangle \right|\right) \right| \nonumber \\
&&~~~\times \left( \left| \langle u,V_m\rangle \right| +\left| \langle P_{E}u,V_m\rangle \right|\right) \nonumber\\
&\leq& 2 \left| \left( ~\left| \langle u,V_m\rangle \right| -\left| \langle P_{E}u,V_m\rangle \right|\right)\right|.\label{inin}\end{eqnarray}
%
\noindent From $$  \langle u,V_m\rangle = \langle P_{E}u,V_m\rangle +  \langle P_{ E^{\bot}} u,V_m\rangle,$$ \noindent 
it readily follows that \begin{equation}\label{inin2} \left|~  \left| \langle u,V_m\rangle \right|-\left| \langle P_{E}u,V_m\rangle \right|~ \right| \leq \left|   \langle P_{ E^{\bot}}u ,V_m\rangle \right|.\end{equation}
(\ref{inin}) and (\ref{inin2}) readily yield
\begin{eqnarray}    \left| ~\left| \langle u,V_m\rangle \right|^2 -\left| \langle P_{E}u,V_m\rangle \right|^2 ~ \right|& \leq& 2 \left|   \langle P_{ E^{\bot}}u ,V_m\rangle \right| \nonumber\\
& \leq & 2 \left\| P_{ E^{\bot}}u \right\|_2     \left\|  V_m \right\|_2   \nonumber  \\& \leq & 2 \left\| P_{ E^{\bot}}u \right\|_2 .
\label{proj1}\end{eqnarray}

\noindent Now, we have that
\begin{eqnarray*}\langle P_{E}u,V_m\rangle&=& \langle u,\xi_m\rangle \langle \xi_m,V_m\rangle\\
 &+&\sum_{\begin{array}{ll}n=1\\n \neq m\end{array} }^{k}\langle u,\xi_n\rangle \langle \xi_n,V_m\rangle \\
&=&\langle u,\xi_m\rangle \langle \xi_m,V_m\rangle+\Delta_m(u) \end{eqnarray*}
where\begin{equation}\label{Delta} \left| \Delta_m(u)\right| \leq
\sum_{\begin{array}{ll}n=1\\n \neq m\end{array} }^{k} \left|\langle \xi_n,V_m\rangle  \right|.\end{equation}
\noindent Then,
$$\left| \langle P_{E}u,V_m\rangle \right|^2=\left|\langle u,\xi_m\rangle  \right|^2 \left|\langle \xi_m,V_m\rangle  \right|^2+\nabla_m(u)$$
\noindent with \begin{equation}\label{nabla} \left| \nabla_m(u)  \right| \leq 2 \left| \Delta_m(u)\right| +\left| \Delta_m(u)\right|^2. \end{equation}
Thus,
\begin{eqnarray}\left| \langle P_{E} u,V_m\rangle \right|^2&=&\alpha \left|  \langle u, \xi_m\rangle \right|^2\nonumber \\ & & + \left(\alpha_m-\alpha\right)\left|  \langle u,\xi_m\rangle \right|^2
 \nonumber \\ & &+ \left(\left|\langle \xi_m,V_m\rangle  \right|^2-\alpha_m\right) \left|  \langle u,\xi_m\rangle \right|^2 \\&& +  \nabla_m(u). \label{proj} 
 \end{eqnarray}
Now, using that $$1= \left\| u \right\|_2^2=\sum_{m=1}^{k} \left|  \langle u, \xi_m \rangle \right|^2 +  \left\| P_{ E^{\bot}}u \right\|_2^2,$$
we have 
\begin{eqnarray*}\sum_{m=1}^{k}\left|  \langle u,V_m\rangle \right|^2 -\alpha &=&\sum_{m=1}^{k} \left[ \left|  \langle u,V_m\rangle \right|^2 -\alpha \left|  \langle u, \xi_m \rangle \right|^2 \right]\\
&& -\alpha \left\| P_{ E^{\bot}}u \right\|_2^2\end{eqnarray*}
and then
\begin{eqnarray}\sum_{m=1}^{k}\left|  \langle u,V_m\rangle \right|^2 -\alpha
 &=&\sum_{m=1}^{k} \left[ 
\left|  \langle u,V_m \rangle \right|^2 - \left| \langle P_{E} u,V_m\rangle \right|^2 \right] \nonumber\\
&&
+\sum_{m=1}^{k} \left[ 
 \left| \langle P_{E} u,V_m\rangle \right|^2
 -\alpha \left|  \langle u, \xi_m\rangle \right|^2 \right] \nonumber\\
&& -\alpha \left\| P_{ E^{\bot}}u \right\|_2^2 \label{de}
.\end{eqnarray}
Using (\ref{de}),
 (\ref{proj1}), (\ref{proj}), (\ref{Delta}) and  (\ref{nabla}), we immediately get that 
\begin{eqnarray}\left| \sum_{m=1}^{k}\left|  \langle u,V_m\rangle \right|^2 -\alpha \right| &\leq &
2 k  \left\| P_{ E^{\bot}} u \right\|_2  \nonumber\\ 
&&+\alpha\left\| P_{ E^{\bot}}u \right\|_2^2 \nonumber\\
&&+\max_{m=1}^{k} \left| \alpha_m-\alpha \right| \nonumber\\
&&+ a_{N},\label{estimgen}\end{eqnarray}
\noindent with\\
~~

 $a_N =\sum_{m=1}^{k}   \left|  \left|\langle \xi_m,V_m\rangle  \right|^2-\alpha_m
 \right| $  $$ + 2 \sum_{\begin{array}{ll}m,n=1\\n \neq m\end{array} }^{k} \left|\langle \xi_n,V_m\rangle  \right|+
\sum_{m=1}^{k}\left(\sum_{\begin{array}{ll}n=1\\n \neq m\end{array} }^{k} \left|\langle \xi_n,V_m\rangle  \right|\right)^2.$$

\noindent Now, if  (\ref{hypconv}) is satisfied, that is,  for each $m,n$ in $\{1,\ldots,k\}^2$,
when $N$ goes to infinity,
$$\left| \langle V_m,\xi_n\rangle \right|^2 {\rightarrow}  \delta_{m,n}\alpha_m,$$ it is clear that 
$a_N$ converges  towards  zero when $N$ goes to infinity.
Lemma \ref{Proj} follows. $\Box$

\begin{lemme}\label {evProj}
Let $ M_N $ be an Hermitian $N \times N$ matrix. Assume that there is a sequence $(n(N))_{N\geq 0}$ in $\mathbb{N}$ and a fixed positive integer number  $k$,
such that for any $l=1,\ldots,k$, $\lambda_{n(N)+l}(M_N)$ converges towards $\rho\in \mathbb{R}$ when $N$ goes to infinity and 
there exists  $ \delta_0>0$ such that,
 for all large N, \begin{equation}\label{deviate}\lambda_{n(N)}(M_N) > \rho + \delta_0 \, \mbox{~and~} \, \lambda_{n(N)+k+1}(M_N) < \rho - \delta_0 
\end{equation}
 (with the conventions that $\lambda_0(M_N)=+\infty$ and $\lambda_{N+1}(M_N)=-\infty$).
For any    $0<\epsilon<\epsilon_0$,  let $M_N(\epsilon)$ be an Hermitian $N \times N$ matrix. Assume that  
 there exists $f_\epsilon \geq 0$, independent of $N$, decreasing to zero when $\epsilon$ decreases to zero,
 such that  for all large  $N$, for any  $0<\epsilon<\epsilon_0$,
\begin{equation}\label{nvelledif}\Vert M_N (\epsilon) - M\Vert \leq f_\epsilon.\end{equation}
Let  $0<\tilde{\epsilon}_0<\epsilon_0$ be such that $f_{\tilde{\epsilon}_0}<\frac{\delta_0}{4}$. 
Then for all large $N$ , for any $0<\epsilon< \tilde{\epsilon}_0$,  for any $l=1,\ldots,k$,  
\begin{equation}\label{convergence}\left|  \lambda_{n(N)+l}(M_N)(\epsilon) -\rho\right|  \leq \frac{\delta_0}{2},\end{equation}
\begin{equation}\label{sep}\lambda_{n(N)}(M_N(\epsilon)) > \rho+ \frac{\delta_0}{2}\mbox{~~and~~}
 \lambda_{n(N)+k +1}(M_N(\epsilon)) < \rho  - \frac{\delta_0}{2},
\end{equation}
\noindent and   for any normalized eigenvector $\xi$ of $M_N$ relative to  the eigenvalue $\lambda_{n(N)+l}(M_N)$ for some  $l$ in $\{1,\ldots,k\}$, 
$$\left\| P_{ E(\epsilon)^{\bot}}\xi\right\|_2 \leq \frac{2}{\delta_0} \left\{ f_\epsilon + \left|  \lambda_{n(N)+l}(M_N)-\rho
\right|\right\},$$
\noindent where  $E(\epsilon)$ denotes the vector subspace generated by the eigenvectors relative to the eigenvalues $\lambda_{n(N)+q}(M_N(\epsilon))$, $q=1,\ldots,k$. 
\end{lemme}
\noindent {\bf Proof:}
\noindent According to Weyl inequalities (see Lemma \ref{Weyl} in the Appendix) and (\ref{nvelledif}),
  for all large $N$, for all   $0<\epsilon<\epsilon_0$,
$$\lambda_{n(N)+l}(M_N(\epsilon))\leq \lambda_{n(N)+l}(M_N) +f_\epsilon,\mbox{~~for~~} l=1, \ldots, k+1,$$
$$\lambda_{n(N)+l}(M_N(\epsilon))\geq \lambda_{n(N)+l}(M_N) -f_\epsilon, \mbox{~~for~~} l=0, \ldots, k.$$

\noindent 
By assumptions of the lemma, 
for all large $N$,  for any $l=1,\ldots,k$,
$$\left|  \lambda_{n(N)+l}(M_N) -\rho\right|  \leq \frac{\delta_0}{4},$$
and
\noindent  for all large N, 
$$\lambda_{n(N)}(M_N) > \rho  + \delta_0 \, \mbox{~and~} \, \lambda_{n(N)+k +1}(M_N) < \rho  - \delta_0 .$$
Hence, choosing  $0<\tilde{\epsilon}_0<\epsilon_0$  such that $f_{\tilde{\epsilon}_0}<\frac{\delta_0}{4}$,  (\ref{sep}) and (\ref{convergence})  readily follow.

 For all large $N$ and any   $0<\epsilon<\tilde{\epsilon}_0$,
let $\xi_1(\epsilon),\ldots, \xi_{k}(\epsilon)$ be an orthonormal basis of $ E(\epsilon)$ such that there exists  an $N \times N$ unitary matrix $V(\epsilon)$ whose $k$ first columns are $\xi_1(\epsilon),\ldots, \xi_{k}(\epsilon)$ and a $(N -k)\times (N-k)$ diagonal matrix $\Lambda_2(\epsilon)$  such that $$M_N(\epsilon)=V(\epsilon)\begin{pmatrix} \Lambda_1(\epsilon) & (0)\\(0)&\Lambda_2(\epsilon)\end{pmatrix}V(\epsilon)^*$$
where $$\Lambda_1(\epsilon)=\mbox{diag}\left(\lambda_{n(N)+1}(M_N(\epsilon)),  \ldots, \lambda_{n(N)+k}(M_N(\epsilon))\right).$$
\noindent Let $\xi$ be a normalized eigenvector  of $M_N$ relative to  the eigenvalue $\lambda_{n(N)+l}(M_N)$ for some  $l$ in $\{1,\ldots,k\}$.
Let us set $$R(\epsilon):=\left( M_N(\epsilon)- \rho I_N\right)\xi.$$
We have 
$$R(\epsilon)=V(\epsilon) \begin{pmatrix}\Lambda_1(\epsilon)- \rho I_{k}& (0)\\(0)&\Lambda_2(\epsilon)- \rho I_{N-k}\end{pmatrix} V(\epsilon)^*\xi.$$
Define the vector $v_1(\epsilon)$ in $\mathbb{C}^{k}$ and the vector $v_2(\epsilon)$ in $\mathbb{C}^{N-k}$
by setting
$$V(\epsilon)^*\xi=\left( \begin{array}{ll} v_1(\epsilon)\\ v_2(\epsilon)\end{array}\right).$$
Note that 
\begin{equation}\label{identite}\left\| v_2(\epsilon)\right\|_2=\left\| P_{ E(\epsilon)^{\bot}} \xi\right\|_2.\end{equation}

\noindent According to (\ref{sep}), for all large $N$,  for any $0<\epsilon < \tilde{\epsilon}_0$,
for all $i\notin \{n(N)+1, \ldots, n(N)+k\}$, $\left| \lambda_i(M_N(\epsilon))- \rho\right| >
\frac{\delta_0}{2} >0$ and therefore
 $\rho$ is not an eigenvalue of $\Lambda_2(\epsilon)$.
  Then  
$$ v_2(\epsilon)= \left(\Lambda_2(\epsilon)- \rho I_{N-k}\right)^{-1} [V^*(\epsilon)R(\epsilon)]_{(N-k)\times 1}$$
\noindent where $[V^*(\epsilon)R(\epsilon)]_{(N-k)\times 1}$ denotes the vector obtained 
from $V^*(\epsilon)R(\epsilon)$ after removing the first $k$ components.
Hence  \begin{equation}\label{v2}\left\| v_2(\epsilon) \right\|_2 \leq \frac{2}{\delta_0} \left\| R(\epsilon) \right\|_2 .\end{equation}
Now, we have
\begin{eqnarray}
\left\| R(\epsilon)\right\|_2 &= &\left\| \left( M_N(\epsilon)-M_N
+M_N -\rho I_N\right)\xi\right\|_2 \nonumber\\
&\leq & \left\| M_N(\epsilon)-M_N \right\|+ \left| \lambda_{n(N)+l}(M_N)-\rho\right| \label{normeR}
\end{eqnarray}

\noindent Lemma \ref{evProj} readily  follows from  (\ref{identite}),  (\ref{v2}),  (\ref{normeR}) and (\ref{nvelledif}). $\Box$\\
~~

Let us define the continuous function $\tau$ on ${\cal O}^{(a)}$ respectively 
${\cal O}^{(m)}\cap ]0;+\infty[$
 by setting
 $$\tau(x)=\left\{\begin{array}{ll} H'(x)  \mbox{~if~} M_N=M_N^W,\\
- \frac{{\cal Z}'(\frac{1}{x})}{x {\cal Z}(\frac{1}{x}) } \mbox{~if~} M_N=M_N^S.
 \end{array}\right.$$
Assume that  $\theta_j $   in $ \Theta_o$ is such that $k_j\neq 1$.
Let us denote by   $V_1(j),\ldots, V_{k_j}(j)$, an orthonormal system of eigenvectors of $A_N$ associated with $\theta_j$. 
 There exists an $N \times N$ unitary matrix $U$ whose $k_j$ first columns are $ V_1(j),\ldots, V_{k_j}(j)$ and a $(N -k_j)\times (N-k_j)$ Hermitian matrix $D$  such that $$A_N=U \begin{pmatrix}\theta_j I_{k_j}& (0) \\ (0)& D \end{pmatrix} U^*.$$
Let us fix  $\epsilon_0$ such that $0< \epsilon_0 < \frac{1}{k_j}  \min_{s=1}^{J}dist(\theta_s, \mbox{supp~}\nu \cup_{i\neq s}\theta_i )$ and $[\theta_j; \theta_j +k_j\epsilon_0]\subset {\cal O}$. For any $0<\epsilon < \epsilon_0$, 
let us consider $$M_N(\epsilon)=\left\{\begin{array}{ll}X_N + A_N(\epsilon)\mbox{~if~} X_N=X_N^W,\\
 A_N(\epsilon)^{\frac{1}{2}}X_N A_N(\epsilon)^{\frac{1}{2}} \mbox{~if~} X_N=X_N^S,\end{array}\right.$$
where $$A_N(\epsilon)=U 
\begin{pmatrix}\mbox{diag}(\theta_j + k_j \epsilon,  \ldots, \theta_j+ 2 \epsilon, \theta_j +  \epsilon)& (0)\\(0)&D
\end{pmatrix}U^*.$$
\noindent Of course for any $0<\epsilon < \epsilon_0$, the limiting spectral distribution of $A_N(\epsilon)$, when $N$ goes to infinity, is the same as the limiting spectral distribution of $A_N$.
Moreover, for all large $N$, the descending ranks $n_{j-1}+1, \ldots , n_{j-1}+k_j$ of $\theta _j$ among the eigenvalues of $A_N$
are the ranks of $\theta_j + k_j \epsilon,  \ldots, \theta_j+ 2 \epsilon,\theta_j + \epsilon$ among the eigenvalues of $A_N(\epsilon)$.
For each $m$ in $\{1,\ldots,k_j\}$, $V_m(j)$ is an eigenvector of $A_N(\epsilon)$ associated with the
 eigenvalue $\theta_j+ (k_j- m+1)\epsilon$ which is of multiplicity one. Note that, since $A_N$ satisfies Assumption A, there exists a constant $C^{'}$ such that  for any $0\leq \epsilon < \epsilon_0$, 
$sup_N \left\| A_N(\epsilon) \right\| \leq C^{'}$. It is easy to see that 
$$\left\| M_N(\epsilon)-M_N\right\| \leq \left\{\begin{array}{ll}  \epsilon ~ k_j   \mbox{~if~} X_N=X_N^W,\\
(C^{'})^{\frac{1}{2}}\frac{k_j}{\sqrt{\theta_j}} \left\| X_N \right\|~ \epsilon \mbox{~if~} X_N=X_N^S.\end{array}\right.$$
According to Theorem 5.11 in \cite{BaiSil06}, $\left\| X_N^S\right\| =c(1+ \frac{1}{\sqrt{c}})^2+o_{a.s,N}(1)$.
Thus in both cases, there exists some constant $C_1$ such that a.s for all large N, 
for any $0< \epsilon< \epsilon_0$,
\begin{equation}\label{normediff} \left\| M_N(\epsilon)-M_N\right\| \leq C_1  \epsilon.\end{equation}

 Let $\xi$ be a normalized  eigenvector of $M_N$ relative to $\lambda_{n_{j-1}+q}(M_N)$ for some $q$ in $\{1,\ldots, k_j\}$ . Let
$0<\tilde{\epsilon}_0<\epsilon_0$ be chosen  such that $C_1 {\tilde{\epsilon}_0}<\frac{\delta_0}{4}$ where $\delta_0$ is defined in Theorem \ref{ThmASCV}.
Using Theorem \ref{ThmASCV} and  (\ref{normediff}), according to Lemma \ref{evProj}, almost surely for all large $N$,
 for any  $0<\epsilon < \tilde{\epsilon}_0$,  the set $\{ 
\lambda_{n_{j-1}+n}(M_N(\epsilon)), n\in \{1,\ldots, k_j\}\}$ is distinct from the set $\{ \lambda_i(M_N(\epsilon)), i\notin \{n_{j-1}+1, \ldots, n_{j-1}+k_j\}\}$ and 
\begin{equation}\label{least2}\left\| P_{ E(\epsilon)^{\bot}}\xi\right\|_2 \leq \frac{2}{\delta_0} \left\{C_1 \epsilon + \left|  \lambda_{n_{j-1}+q}(M_N)-\rho_{\theta_j}
\right|\right\},\end{equation}
where  $E(\epsilon)$ denotes the vector subspace generated by the eigenvectors relative to the eigenvalues $\lambda_{n_{j-1}+n}(M_N(\epsilon))$, $n=1,\ldots,k_j$. 
Define \begin{equation}\label{eta}\iota(\epsilon)= \frac{2 (2k_j+\tau(\theta_j))}{\delta_0} C_1 \epsilon + \max_{m=1}^{k_j} \left| \tau(\theta_j + m\epsilon)-\tau(\theta_j)\right|.\end{equation}
\noindent For any $\zeta >0$, choose and fix $0<\epsilon= \epsilon_1<\tilde{\epsilon}_0$ such that \begin{equation} \label{iota}0<\iota(\epsilon_1)< \frac{\zeta}{2}\end{equation}
(using the continuity of the function $\tau$ at the point $\theta_j$). \\
~~

\noindent Now, each $\theta_j+l \epsilon_1$, $l\in \{1,\ldots, k_j\}$, is a spike of $A_N(\epsilon_1)$ with  multiplicity one.
According to Theorem \ref{ThmASCV},  for  $n\in \{1,\ldots, k_j\}$,
 $\lambda_{n_{j-1}+n}(M_N(\epsilon_1))$ asymptotically separates from the rest of the spectrum and converges almost surely towards 
$H(\theta_j +(k_j-n+1) \epsilon_1)$; moreover, if $\xi_n(\epsilon_1,j)$ denotes a normalized eigenvector associated to $\lambda_{n_{j-1}+n}(M_N(\epsilon_1))$,
 Proposition \ref{cvev}
implies that
$$\left| \langle V_m(j),\xi_n(\epsilon_1,j)\rangle \right|^2 \stackrel{a.s}{\rightarrow}  \delta_{m,n}\tau(\theta_j+(k_j-m+1)\epsilon_1).$$
According to Lemma \ref{Proj},  there exists  a random variable $a_N(\epsilon_1) \geq 0$, converging almost surely  to zero when $N$ goes to infinity,
 such that, almost surely, for all large $N$,
\begin{eqnarray*}\left|  \left\| P_{\mbox{Ker~}(\theta_j I_N-A_N)} \xi \right\|_2^2-\tau(\theta_j)\right| \leq& (2k_j+\tau(\theta_j))  \left\|P_{E(\epsilon_1)^{\bot}} \xi \right\|_2 +a_N(\epsilon_1)\nonumber \\
&+\max_{m=1}^{k_j} \left| \tau(\theta_j + m\epsilon_1)-\tau(\theta_j)\right|.\end{eqnarray*}
The last inequality, (\ref{least2}) and (\ref{eta}) readily yield  that, almost surely, for all large $N$,

\begin{eqnarray*}\label{least}\left|  \left\| P_{\mbox{Ker~}(\theta_j I_N-A_N)} \xi \right\|_2^2-\tau(\theta_j)\right| &\leq&  \frac{2}{\delta_0}(2k_j+\tau(\theta_j))   \left|  \lambda_{n_{j-1}+q}(M_N)-\rho_{\theta_j}
\right| \\ && +a_N(\epsilon_1) +\iota(\epsilon_1).\end{eqnarray*}

\noindent Therefore, (\ref{iota}), the  almost surely  convergence of  $ \lambda_{n_{j-1}+q}(M_N)$ towards $\rho_{\theta_j}$ and of  $a_N(\epsilon_1)$
towards zero imply that,
almost surely, for all large $N$, 
$$\left|\left\| P_{\mbox{Ker~}(\theta_j I_N-A_N)} \xi \right\|^2_2 -\tau(\theta_j) \right| 
\leq 
\zeta $$
\noindent and the proof is  complete.\hspace*{8cm} $\Box$

\section{Proof of Proposition \ref{cvev}}
Since the proof of Proposition \ref{cvev} is exactly the same for the deformed Wigner model and the sample covariance matrix, in order to present a unified approach  of this proof, we adopt in this section the notations of the first column of the  array (\ref{correspondance}) standing for both the corresponding elements of the second column (deformed Wigner matrix case) and the third column (sample covariance matrix case). We postpone in a later subsection the technical results that need a specific study for each model in order to not lose the thread of this common proof.

\subsection{Restriction  to the asymptotic behavior of some expectation $\E\Tr \left[h(M_N) f(A_N)\right]$}
The aim of this first step of the proof of Proposition \ref{cvev} is to reduce the study of the asymptotic behaviour of  $\sum_{n=1}^{k_j}\left\| P_{\mbox{Ker~}(\theta_l I_N-A_N)} \xi_n(j) \right\|^2_2$ to the one of the expectation $\E\Tr \left[h(M_N) f(A_N)\right]$ for some  functions $f$ and $h$ respectively concentrated on a neighborhood of $\theta_l$ and $\rho_{\theta_j}$. We will use the convergence results on eigenvalues described in Theorem \ref{cvev} above and concentration inequalities presented in the Appendix.\\
P. Biane  already suggested in \cite{PB} to evaluate the moduli of the Hermitian inner products of the eigenvectors on test functions. 
Indeed, for any smooth function $h$ and $f$ on $\mathbb{R}$, denoted by $u_1,\ldots, u_N$ (resp. $w_1,\ldots, w_N$), the eigenvectors associated with $\lambda_1(A_N), \ldots, \lambda_N(A_N)$ (resp. $\lambda_1(M_N), \ldots, \lambda_N(M_N)$),  one can easily check that
$$\Tr \left[h(M_N) f(A_N)\right] =\sum_{k,i} h(\lambda_k(M_N)) f(\lambda_i(A_N)) \vert \langle u_i,w_k \rangle \vert^2.$$
Thus, since $\theta_l$ on one hand  and the $\lambda_{n_{j-1}+n}(M_N)$, $n=1,\ldots,k_j$, on the other hand, asymptotically separate from the rest of the spectrum of respectively $A_N$ and $M_N$, 
a fit choice of $h$ and $f$ will allow the  study of the restrictive sum $\sum_{n=1}^{k_j}\left\| P_{\mbox{Ker~}(\theta_l I_N-A_N)} \xi_n(j) \right\|^2_2$.\\

\noindent Let us fix $$0< \eta < \frac{1}{2} \min_{s=1}^{J}dist(\theta_s, \mbox{supp~}\nu \cup_{i\neq s}\theta_i )$$ and
for any $l=1,\ldots,J$,
choose $f_{\eta,l}$ in ${\cal C}^\infty (\mathbb{R}, \mathbb{R})$ with support in 
$[\theta_l -\eta,\theta_l+\eta]$  such that $f_{\eta,l}(\theta_l)=1$ and $0 \leq f_{\eta,l} \leq 1$.\\
~\\
For any $\theta_i \in {\cal O}^{(a)}$, according to Remark \ref{oainclu}, $\theta_i$ belongs to $\overline{\Omega_{\sigma,\nu}}$ and according to Proposition \ref{homeo},  $F_{\sigma,\nu}^{(a)}(\rho_{\theta_i})=\theta_i$.
For any $\theta_i \in {\cal O}^{(m)}\cap ]0;+\infty[$, according to Remark \ref{nonnul}, $\rho_{\theta_i}= {\cal Z}(\frac{1}{\theta_i})>0$
and according to Remark \ref{CScomp}, $g_{\tau_{c,\nu}}(\rho_{\theta_i}) =\frac{1}{\theta_i}$ so that 
$\frac{1}{F_{c,\nu}^{(m)}(\frac{1}{\rho_{\theta_i}})} =\theta_i$.\\
According to  Theorem \ref{ThmASCV}, there exists 
$\delta_0>0$  such that almost surely for all large N, for all $\theta_j$ in $\Theta_o$, 
$$\lambda_{n_{j-1}}(M_N) > \rho _{\theta _j} + \delta_0 \, \mbox{~and~} \, \lambda_{n_{j-1}+k_j +1}(M_N) < \rho _{\theta _j} - \delta_0 .$$ 

\noindent Let us fix $$0< \delta < \frac{1}{3}\min \{\delta_0, \mbox{~dist~}(\rho_{\theta_s}, {\cal S}), \vert \rho_{\theta_s}- \rho_{\theta_{t}}\vert, s\neq t,(\theta_s,\theta_t) \in \Theta_o^2\},$$
where $${\cal S}=   \left\{ \begin{array}{ll} \mbox{supp~}(\mu\boxplus \nu) \mbox{~if ~} M_N =M_N^W\\
\mbox{supp~}(\tau_{c,\nu})\cup \{0\}\mbox{~if ~} M_N =M_N^S \end{array} \right. .$$
\noindent For any $j$ such that $\theta_j \in \Theta_o$, choose $h_{\delta,j}$ in ${\cal C}^\infty (\mathbb{R}, \mathbb{R})$ with support in 
$[\rho_{\theta_j} -\delta,\rho_{\theta_j}+\delta]$  such that $h_{\delta,j} \equiv 1$ on $[\rho_{\theta_j} -\frac{\delta}{2},\rho_{\theta_j}+\frac{\delta}{2}]$ and $0 \leq h_{\delta,j}\leq 1$. 

\begin{lemme} \label{vecttrace}
When $N$ goes to infinity $$Tr \left[h_{\delta,j}(M_N) f_{\eta,l}(A_N)\right]-\sum_{n=1}^{k_j}\left\| P_{\mbox{Ker~}(\theta_l I_N-A_N)} \xi_n(j) \right\|^2_2 \stackrel{a.s}{\rightarrow}0.$$
\end{lemme}
\noindent {\bf Proof:} According to Theorem \ref{ThmASCV},
there exists some set $\Omega$ of probability one 
 such that on $\Omega$,  for all large $N$, $\forall i=1,\ldots,k_j,$
$$
\vert \lambda_{n_{j-1}+i}(M_N) -\rho _{\theta _j}\vert < \frac{\delta}{2} $$
 $$ \lambda_{n_{j-1}}(M_N)\geq \rho_{\theta_j} +  \delta, \, \lambda_{n_{j-1}+ k_j+1}(M_N) \leq \rho_{\theta_j} -  \delta .$$
Using also the assumption (\ref{univconv}) on the $\beta_i(N)$'s, we have that  on $\Omega$, for all large $N$, $$\sum_{n=1}^{k_j}\left\| P_{\mbox{Ker~}(\theta_l I_N-A_N)} \xi_n(j) \right\|^2_2 = Tr \left[h_{\delta,j}(M_N) f_{\eta,l}(A_N)\right].$$
Hence,
Lemma \ref{vecttrace} follows. $\Box$\\

Now, according to Lemma \ref{Herbst}, Remark \ref{multiple} and Lemma \ref{Lipschitz}, the random variables $F^W_N=\Tr \left[ h_{\delta,j} (M_N^W) f_{\eta,l}(A_N)\right]$
and $F^S_N=\Tr \left[ h_{\delta,j} (M_N^S) f_{\eta,l}(A_N)\right]$ 
satisfy respectively  the following concentration inequalities
$$\forall \epsilon> 0, \, \mathbb{P}\left( \vert F^W_N-\E(F^W_N) \vert > \epsilon \right) \leq K_1 \exp\left(-\frac{\epsilon\sqrt{ N}}{K_2 \sqrt{ C_{PI} k_l } \Vert h_{\delta,j} \Vert_{Lip}}\right).$$
$$\forall \epsilon> 0, \, \mathbb{P}\left( \vert F^S_N-\E(F^S_N) \vert > \epsilon \right) \leq K_1 \exp\left(-\frac{\epsilon \sqrt{ p}}{K_2 \sqrt{C C_{PI} k_l } \Vert \tilde{h}_{\delta,j} \Vert_{Lip}}\right).$$
\noindent (The constants have been introduced in  Lemma \ref{Herbst} and Lemma \ref{Lipschitz}).
By Borel-Cantelli Lemma, we can readily deduce the following lemma.  
\begin{lemme}\label{conc}
$$\Tr \left[ h_{\delta,j} (M_N) f_{\eta,l}(A_N)\right]- \mathbb{E}\left[ \Tr \left[ h_{\delta,j} (M_N) f_{\eta,l}(A_N)\right] \right]\stackrel{a.s}{\rightarrow}0.$$

\end{lemme}

Lemma \ref{vecttrace} and Lemma \ref{conc} allow us to conclude this first step of the proof by the following result.
\begin{proposition}\label{premetap}
For any $\theta_j$ in $\Theta_o$ and any $\theta_l$ in $\Theta$, when $N$ goes to infinity
$$\sum_{n=1}^{k_j}\left\| P_{\mbox{Ker~}(\theta_l I_N-A_N)} \xi_n(j) \right\|^2_2 -\mathbb{E}\left[ \Tr \left[ h_{\delta,j} (M_N) f_{\eta,l}(A_N)\right] \right]\stackrel{a.s}{\rightarrow}0.$$
\end{proposition}
\subsection{Making use of estimations of the resolvent}
The basic idea of this second step of the proof of Proposition \ref{cvev} is to approximate  the function $h_{\delta,j}$ by its convolution by the Poisson Kernel in order to exhibit the resolvent of the deformed model and then use sharp estimations on this resolvent.
\begin{lemme}\label{approxcauchy}
For any continuous function $h$ with compact support and any  bounded continuous  function $\phi$,
$$\mathbb{E}\left[ \Tr \left[ h (M_N) \phi (A_N)\right] \right]=- \lim_{y\rightarrow 0^{+}}\frac{1}{\pi} \Im \int  \mathbb{E} \left( \Tr \left[G_{N}(t+iy)\phi (A_N)\right] \right)h(t) dt $$
\noindent where $G_N(z) =(zI_N -M_N)^{-1}$.
\end{lemme}
\noindent {\bf Proof}
Let us denote by $P_y$ the Poisson kernel $$P_y(x) =\frac{y}{\pi (x^2 +y ^2)}, \, x\in \mathbb{R}, \, y>0.$$
We have 
\begin{eqnarray*}h(x)&= & \lim_{y\rightarrow 0^{+}} h * P_y(x)\\
&=&\lim_{y\rightarrow 0^{+}}\frac{1}{\pi} \int \frac{y h(t)}{(x-t)^2 +y^2} dt \\
&=& \lim_{y\rightarrow 0^{+}}\frac{1}{\pi} \int \Im\frac{h(t)}{x-iy-t} dt.
\end{eqnarray*}
Thus, for any fixed $N$,  $$h(M_N) = - \lim_{y\rightarrow 0^{+}}\frac{1}{\pi} \int \Im G_{N}(t+iy) h(t) dt, $$
and since $\Vert h * P_y \Vert_\infty \leq \Vert h \Vert_\infty$ we have
$\Vert \int \Im G_{N}(t+iy) h(t) dt\Vert \leq \Vert h \Vert_\infty$.
Then, the result readily follows by dominated convergence Theorem and Fubini's Theorem. $\Box$\\

\noindent 
Let $U$ be a unitary matrix and $$D={\rm diag}(\gamma_1, \ldots ,\gamma_N )$$
such that 
$$A_N= U^*DU$$
Let $G$ stand for $G_N$ and $g$ stand for $g_{\mu_{\mbox{\tiny{LSD}}}}$.
Consider $\tilde{ G }= UGU^*$.
For any continuous function $\phi$, \begin{equation}\label{GtildeA}\Tr \left[G_{M_N}(z)\phi (A_N)\right]= \sum_{k=1}^N \phi (\gamma_k)\tilde{G}_{kk}(z).\end{equation}
The following result is fundamental in our approach.
\begin{proposition}\label{estimfonda}
 There is a polynomial $P$ with nonnegative coefficients,  a sequence $a_N$ of nonnegative
real numbers converging to zero when $N$ goes to infinity 
and some nonnegative integer number $\alpha$, 
such that for any $k$ in $\{1, \ldots,N\}$, for all  $z\in \mathbb{C}\setminus \mathbb{R}$,
\begin{equation} 
\E (\tilde G_{kk}(z)) = \tilde{\Phi}_k(z)
+\Delta_{k,N}(z),
\end{equation}
with $$\left| \Delta_{k,N} (z)\right| \leq (1+\vert z\vert)^\alpha P(\vert \Im z \vert^{-1})a_N,$$
where $$\tilde{\Phi}_k(z)=\left\{ \begin{array}{ll}
\frac{1}{z-\sigma^2 g(z)-\gamma_k}=\frac{1}{F^{(a)}_{\sigma,\nu}(z)-\gamma_k}\, \hspace*{1cm} \mbox{~if~} \, M_N=M_N^W,\\
\frac{1}{z-\gamma_k(1-c+czg(z))}=\frac{1}{z ( 1 -\gamma_k F^{(m)}_{c,\nu}(\frac{1}{z}))}\, \mbox{~if~} \, M_N=M_N^S.\end{array} \right.$$
\end{proposition}
Note that $\tilde{\Phi}_k(z)$ is well defined for any $z\in \mathbb{C}\setminus \mathbb{R}$ since 
\begin{equation}\label{majim}\vert \Im[z-\sigma^2 g(z)-\gamma_k]\vert = \vert \Im(z)\vert  \left[1+ \sigma^2\int \frac{1}{\vert z-x\vert^2}d\mu_{\sigma}\boxplus \nu(x)\right]\geq \vert \Im (z)\vert >0,\end{equation}
\begin{equation}\label{majim2}\vert \Im[z-\gamma_k(1-c+czg(z))]\vert = \vert \Im(z)\vert  \left[1+ c\gamma_k\int \frac{x}{\vert z-x\vert^2}d (\mu_{\mbox{\tiny{MP}}, {c}}\boxtimes \nu)(x)\right]\geq \vert \Im (z)\vert >0.\end{equation}

According to (\ref{estimfondakW}) and Proposition \ref{estimresolvwish},
there exists a polynomial $P$ with nonnegative coefficients and a sequence $a_N$ of nonnegative
real numbers converging to zero when $N$ goes to infinity such that, for any $k$ in $\{1, \ldots,N\}$, for any $z\in \mathbb{C}\setminus \mathbb{R}$ \begin{equation} \label{estimfondaN}
\E (\tilde G_{kk}(z)) = \tilde{\Phi}_{N,k}(z)
+R_{k,N}(z),
\end{equation}
with $$\vert R_{k,N}(z) \vert \leq (1+\vert z\vert)^2 P(\vert \Im z \vert^{-1})a_N,$$
where $$\tilde{\Phi}_{N,k}(z)=\left\{ \begin{array}{ll}
\frac{1}{z-\sigma^2 g_N^W(z)-\gamma_k}, \hspace*{1cm} \mbox{~if~} \, M_N=M_N^W,\\
\frac{1}{z-\gamma_k(1-\frac{N}{p}+\frac{N}{p}zg_N^S(z))}=\, \mbox{~if~} \, M_N=M_N^S.\end{array} \right.$$

In order to deduce Proposition \ref{estimfonda},  we will need the following description of the convergence of $g_N(z)$ towards $g(z)$. 
\begin{proposition}\label{estimdif}
There exists a polynomial $R$ with nonnegative coefficients, a sequence $a_N$ of positive numbers converging towards zero and some nonnegative integer number $\alpha$ such that,  for all $z \in \mathbb{C}\setminus \mathbb{R}$,
 \begin{equation}\label{dif}\vert {g}_{N}(z) -g(z)\vert \leq (1+ \vert z \vert)^\alpha R(\vert \Im z \vert^{-1}) a_N.\end{equation}
\end{proposition}
\noindent{\bf Proof}: 
1) The deformed Wigner model case.\\
Denote by $\tilde{g}_{N}$ the Stieltjes transform of $\mu_\sigma \boxplus \mu_{A_N}$.
Since  we have already proved in \cite{CDFF10} that there exists a polynomial $S$ with nonnegative coefficients such that
for all $z \in \mathbb{C}\setminus \mathbb{R}$, \begin{equation}\vert \tilde{g}_{N}(z) -g_N(z)\vert \leq \frac{S(\vert \Im z \vert^{-1})}{N},\end{equation}
the result will readily follow if we prove that
 there exists a polynomial $T$ with nonnegative coefficients and a sequence $b_N$ of positive numbers converging towards zero such that,  for all $z \in \mathbb{C}^+$,
 \begin{equation}\label{dif2}\vert \tilde{g}_{N}(z) -g(z)\vert \leq T(\vert \Im z \vert^{-1}) b_N.\end{equation}
The proof of (\ref{dif2}) follows the lines of Section 4 in \cite{CDFF10}. 
For a fixed $z\in \C^+$, according to Proposition \ref{homeo}, we have the subordination equations: 
\begin{equation} \label{eqtilde} \tilde g_N(z) = g_{\mu _{A_N}}(F^{(a)}_{\sigma,\mu _{A_N}}(z))=g_{\mu _{A_N}}(z - \sigma ^2\tilde g_N(z)),\end{equation}
\begin{equation} \label{eq}  g(z) = g_{\nu}(F^{(a)}_{\sigma,\nu}(z))=g_{\nu}(z - \sigma ^2 g(z)).\end{equation}
Moreover, using  Lemma \ref{cvAN}  and $\Im (z - \sigma ^2g(z)) \geq  \Im z $ , 
we deduce from (\ref{eq}) that 
\begin{equation} \label{appsub} g(z) = g_{\mu _{A_N}}(z - \sigma ^2g(z)) + \Delta_N(z) \end{equation}
\noindent with 
$$\vert \Delta_N (z)\vert \leq v_N(1) P_1(|\Im z|^{-1}),$$
\noindent where $P_1$ is a polynomial with nonnegative coefficients and  $v_N(1)$ is a sequence of  positive numbers converging towards zero. \\
Since $z - \sigma ^2g(z)\in \C^+$, $z'\in \C$ is 
well-defined by the formula : 
$$z':=H_{\sigma , \mu _{A_N}}(z - \sigma ^2g(z)),$$
where $H_{\sigma , \mu _{A_N}}$ is defined by (\ref{defdeH}) replacing $\nu$ by $\mu_{A_N}$.
One has 
\begin{eqnarray*}
\vert z'-z \vert&=&\vert-\sigma ^2(g(z)-g_{\mu _{A_N}}(z - \sigma ^2g(z)))\vert \\
    &\leq &\sigma^2 v_N(1) P_1(|\Im z|^{-1}).
\end{eqnarray*}

\begin{itemize}
\item
If $$\frac{|\Im z|}{2} \leq \sigma^2 v_N(1) P_1(|\Im z|^{-1}),$$ 
or equivalently 
\begin{equation} \label{1=O(1/N)}
1\leq {2\sigma^2|\Im z|^{-1}P_1(|\Im z|^{-1})}v_N(1),
\end{equation}
\begin{eqnarray*}
|g(z)-\tilde g_N(z)|&\leq &\frac{2}{|\Im z|}\\
&\leq &{4\sigma^2|\Im z|^{-2}P_1(|\Im z|^{-1})}v_N(1).
\end{eqnarray*}
\item If $$\frac{|\Im z|}{2}> \sigma^2 v_N(1) P_1(|\Im z|^{-1}),$$ 
one has : 
$$ |\Im z'-\Im z|\leq |z'-z|\leq \frac{|\Im z|}{2}$$
which implies $\Im z'\geq \frac{\Im z}{2}$ and 
therefore $z'\in \C^+$. 
Hence, according to (\ref{imageinverse}), it follows that
$z-\sigma ^2g(z)\in \Omega_{\sigma , \mu _{A_N}}$ (where $ \Omega_{\sigma , \mu _{A_N}}$ is defined by (\ref{Ome}) replacing
$\nu$ by  $\mu _{A_N}$)
so that $F^{(a)}_{\sigma , \mu _{A_N}}(z')=z-\sigma ^2g(z)$. 

Thus, the approximative equation \eqref{appsub} may be rewritten 
$$g(z) =g_{\mu _{A_N}}(F^{(a)}_{\sigma , \mu _{A_N}}(z')) + \Delta_N(z)$$
and then, using the subordination equation \eqref{eqtilde} 
$$
g(z)=\tilde g_N(z')+
\Delta_N(z).$$

Moreover,
\begin{eqnarray*}\vert \tilde{g}_N(z')-\tilde{g}_N(z)\vert&=& \vert (z-z')\int_{\mathbb{R}}\frac{d(\mu_\sigma \boxplus \mu_{A_N})(x)}{(z'-x)(z-x)}\vert\\
&\leq & 2\sigma^2 v_N(1) \vert \Im z \vert^{-2} P_1(|\Im z|^{-1}).\end{eqnarray*}

Hence 
\begin{eqnarray*}
|g(z)-\tilde g_N(z)|&\leq &|g(z)-\tilde g_N(z')|+|\tilde g_N(z')-\tilde g_N(z)|\\ & \leq & (2\sigma^2  \vert \Im z \vert^{-2}+1)v_N(1) P_1(|\Im z|^{-1})
\end{eqnarray*}
\end{itemize}
Finally we get that for all $z \in \mathbb{C}^+$, 
$$|g(z)-\tilde g_N(z)|\leq  (4\sigma^2  \vert \Im z \vert^{-2}+1)v_N(1) P_1(|\Im z|^{-1})$$
\noindent so that (\ref{dif}) is satisfied in the deformed Wigner model setting with $a_N=v_N(1)$ and $R(x)=(4\sigma^2  x^{2}+1) P_1(x)$, $\alpha=0$.\\

\noindent 2) The sample covariance matrix setting\\
Let $z$ be in $\mathbb{C}^+$.
Note that it is obviously equivalent to prove such an estimation for $\left|\underline{g}_N(z)- g_{\tau_{c,\nu}}(z)\right|$ where $$\underline{g}_N(z) =\frac{1-\frac{N}{p}}{z} + \frac{N}{p} g_N(z)$$
is the expected value of the Stieltjes transform of the spectral measure of $\underline{M}_N^S=\frac{1}{p} B_N^* A_N B_N$.
In the following any $P_i$ will denote a polynomial with nonnegative coefficients and $a_N(i)$ will denote a sequence of nonnegative numbers converging towards zero when $N$ goes to infinity.
Letting the sum running over $k$ in (\ref{submat}) and dividing by $N$ we have
$${g}_N(z) = \int \frac{d\mu_{A_N}(t)}{z(1- t \underline{g}_N(z))} + \Delta_N(z)$$
\noindent where
$$\vert \Delta_N(z) \vert \leq (1+\vert z\vert)^2 P_1(\vert \Im z \vert^{-1})a_N(1).$$
It readily follows that 
$$
\underline{g}_N(z) = \frac{\underline{g}_N(z)}{z} {\cal Z}_{\frac{N}{p},\mu_{A_N}}(\underline{g}_N(z)) + \frac{N}{p}\Delta_N(z)$$
\noindent where ${\cal Z}_{\frac{N}{p},\mu_{A_N}}$ is defined by (\ref{defzcnu}) replacing $\nu$ by $\mu_{A_N}$ and $c$ by 
$\frac{N}{p}$.
\noindent Then using Lemma \ref{MajinversegN} we deduce that
\begin{equation}\label{zzprime}{\cal Z}_{ \frac{N}{p},\mu_{A_N}}(\underline{g}_N(z)) =z+ R_N(z)\end{equation}
\noindent where for all large $N$ $$\vert R_N(z) \vert \leq (1+\vert z\vert)^5 P_2(\vert \Im z \vert^{-1})a_N(1)$$
\noindent  On the other hand, using (\ref{hypA2}) and Lemma \ref{MajinversegN} we have \begin{equation}\label{zzprime2} {\cal Z}_{\frac{N}{p},\mu_{A_N}}(\underline{g}_N(z)) ={\cal Z}(\underline{g}_N(z)) + Q_N(z)\end{equation}
where $$\vert Q_N(z) \vert \leq (1+\vert z\vert)^q P_3(\vert \Im z \vert^{-1})a_N(2)$$ for some nonnegative integer number $q$.
We readily deduce from (\ref{zzprime}) and (\ref{zzprime2}) that 
$${\cal Z}(\underline{g}_N(z)) =z + T_N(z)$$
where $$\vert T_N(z) \vert \leq (1+\vert z\vert)^{\alpha} P_4(\vert \Im z \vert^{-1})a_N(3)$$ for some nonnegative integer number $\alpha$.

\noindent Set $$z^{'}={\cal Z}(\underline{g}_N(z)).$$
\begin{itemize}
\item
If $$\frac{|\Im z|}{2} \leq (1+\vert z\vert)^{\alpha} P_4(\vert \Im z \vert^{-1})a_N(3),$$ 
or equivalently 
$$ 1\leq 2|\Im z|^{-1}(1+\vert z\vert)^{\alpha} P_4(\vert \Im z \vert^{-1})a_N(3),$$
\noindent then
\begin{eqnarray*}
|g_{\tau_{c,\nu}}(z)-\underline{g}_N(z)\vert &\leq &\frac{2}{|\Im z|}\\
&\leq &4|\Im z|^{-2}(1+\vert z\vert)^{\alpha} P_4(\vert \Im z \vert^{-1})a_N(3).
\end{eqnarray*}
\item If $$\frac{|\Im z|}{2}> (1+\vert z\vert)^{\alpha} P_4(\vert \Im z \vert^{-1})a_N(3),$$ 
one has : 
$$ |\Im z'-\Im z|\leq |z'-z|< \frac{|\Im z|}{2}$$
which implies $\Im z'\geq \frac{\Im z}{2}$ and 
therefore $z'\in \mathbb{C}^+$. 
Note that $\underline{g}_N(z)$ satisfied the equation 
$$Z =\frac{1}{z^{'}-c\int \frac{t}{1-tZ}d\nu(t)}$$ and since $g_{\tau_{c,\nu}}(z^{'})$ is the unique solution in 
$\mathbb{C}^{-}$ of the latter equation, we can deduce that 
$g_{\tau_{c,\nu}}(z^{'})=\underline{g}_N(z)$.

Hence 
\begin{eqnarray*}\left|\underline{g}_N(z)- g_{\tau_{c,\nu}}(z)\right|&=& \left|g_{\tau_{c,\nu}}(z^{'})-g_{\tau_{c,\nu}}(z)\right|\\
&\leq & \vert \Im z \vert^{-1} \vert \Im z'\vert^{-1} \vert z - z'\vert \\
&\leq& 2\vert \Im z \vert^{-2}(1+\vert z\vert)^{\alpha} P_4(\vert \Im z \vert^{-1})a_N(3).
\end{eqnarray*}

\end{itemize}
Finally we get that for all $z \in \mathbb{C}^+$, 
$$\left|\underline{g}_N(z)- g_{\tau_{c,\nu}}(z)\right|\leq  4\vert \Im z \vert^{-2}(1+\vert z\vert)^{\alpha} P_4(\vert \Im z \vert^{-1})a_N(3)$$
\noindent so that (\ref{dif}) is satisfied in the sample covariance matrix setting .\\

$\Box$\\

Now, Proposition \ref{estimfonda} readily follows from (\ref{estimfondaN}) and Proposition \ref{estimdif}
using (\ref{image}), (\ref{majim2}) and (\ref{imageW}) and (\ref{majim}). $\Box$\\

Thus, for any $l=1,\ldots,J$, (\ref{GtildeA}) and Proposition \ref{estimfonda} yield that  
for all large $N$, for all  $z \in \mathbb{C} \setminus  \mathbb{R}$, \begin{equation}\label{ecrit} \mathbb{E} \left( \Tr \left[G_{M_N}(z)f_{\eta,l} (A_N)\right] \right)
= \phi_l (z) + \Delta_N(z)\end{equation}
where $\phi_l$ is the following analytic function on $\mathbb{C} \setminus  \mathbb{R}$:
$${\phi}_l(z)=\left\{ \begin{array}{ll}
\frac{k_l}{z-\sigma^2 g(z)-\theta_l}=\frac{k_l}{F^{(a)}_{\sigma,\nu}(z)-\theta_l}\, \hspace*{1cm} \mbox{~if~} \, M_N=M_N^W,\\
\frac{k_l}{z-\theta_l(1-c+czg(z))}=\frac{k_l}{z ( 1 -\theta_l F^{(m)}_{c,\nu}(\frac{1}{z}))}\, \mbox{~if~} \, M_N=M_N^S,\end{array} \right.$$
\noindent and 
$$ \vert \Delta_N(z) \vert \leq a_N (1 + \vert z \vert )^\alpha P(\frac{1}{\vert \Im z \vert})$$
\noindent for some polynomial $P$ with nonnegative coefficients, some sequence $a_N$ of positive
real numbers converging to zero when $N$ goes to infinity and some nonnegative integer number $\alpha$.
 
\noindent According to Lemma \ref{HT}, we  have
$$\limsup_{y \rightarrow 0^+} ~ (a_N)^{-1}\vert  \int h_{\delta,j}(t) \Delta_N(t+iy)dt \vert <+\infty$$
so that \begin{equation}\label{reste}\lim_{N\rightarrow + \infty}\limsup_{y \rightarrow 0^+} \vert \int h_{\delta,j}(t) \Delta_N(t+iy)dt \vert =0.\end{equation}
\begin{lemme}\label{pourholomorphe}
Let $\theta_j $ be in $ \Theta_o$.
\begin{itemize}
\item[(1)] Set $\rho_{\theta_j}=H(\theta_j)$.
If $\theta_l \in \Theta\setminus\{\theta_j\}$, the map
$(x,y) \mapsto F^{(a)}_{\sigma,\nu}(x+iy) -\theta_l$  does not vanish on $[\rho_{\theta_j}-2\delta; \rho_{\theta_j}+2\delta]\times \mathbb{R}$.\\
The only vanishing point in $[\rho_{\theta_j}-2\delta; \rho_{\theta_j}+2\delta]\times \mathbb{R}$ of the map $(x,y) \mapsto F^{(a)}_{\sigma,\nu}(x+iy) -\theta_j$  is $(\rho_{\theta_j},0)$.
\item[(2)]   Set $\rho_{\theta_j}={\cal Z}(\frac{1}{\theta_j})$. If $\theta_l \in \Theta\setminus\{\theta_j\}$, the map $(x,y) \mapsto (x+iy) (1- \theta_l g_{\tau_{c,\nu}}(x+iy))$ does not vanish on $[\rho_{\theta_j}-2\delta; \rho_{\theta_j}+2\delta]\times \mathbb{R}$.\\
The only vanishing point in $[\rho_{\theta_j}-2\delta; \rho_{\theta_j}+2\delta]\times \mathbb{R}$ of the map $(x,y) \mapsto (x+iy) (1- \theta_j g_{\tau_{c,\nu}}(x+iy))$ is $(\rho_{\theta_j},0)$.

\end{itemize}
\end{lemme}

\noindent {\bf Proof:}
Note  that if  $y\neq 0$,  for any $x$, the imaginary part of $F^{(a)}_{\sigma,\nu}(x+iy) -\theta_l$ and $1- \theta_l g_{\tau_{c,\nu}}(x+iy)$ is 
nonnull
so that we will focus on the case $y=0$.\\
\noindent{ Proof of (1):}
\begin{itemize}
\item Assume $\theta_l \notin \Theta_o$.
First, if $H'(\theta_l)<0$, according to (\ref{imagereelle}), 
$\theta_l $ does not belong to $F^{(a)}_{\sigma,\nu}(\mathbb{R})$  so that the conclusion of Lemma \ref{pourholomorphe} (1) is  true.\\
Now assume that $H'(\theta_l)=0$. According to (\ref{imagereelle}),  $\theta_l\in \partial \Omega_{\sigma,\nu}=F^{(a)}_{\sigma,\nu}(\mathbb{R})$,
and,
by Proposition \ref{homeo}, $F^{(a)}_{\sigma,\nu}(x) -\theta_l=0$ implies $x=H(\theta_l)$.
 For any $u \in {\cal O}^{(a)}$,
we have $H(\theta_l)\neq H(u)$. Indeed, for any $u \in  {\cal O}^{(a)}$, there exists $u_1, u_2$ such that $[u_1;u_2] \subset  {\cal O}^{(a)}$ and $u_1< u < u_2$. Since $H$ is globally nondecreasing on 
$\{v \in ^c \mbox{supp~}(\nu), \, H'(v)\geq 0\}$ 
(see Remark \ref{Hcroit}), we have if $\theta_l <u$, 
$H(\theta_l) \leq H(u_1) < H(u)$ 
and if  $\theta_l >u$, $H(u) < H(u_2) \leq H(\theta_l)$.
It follows, according to  (\ref{car1}),  that $H(\theta_l)$ belongs to supp$(\mu\boxplus \sigma)$  and therefore cannot belong to $[\rho_{\theta_j}-2\delta; \rho_{\theta_j}+2\delta]$ so that the conclusion of Lemma \ref{pourholomorphe} (1) is true. 

\item 
Let us consider now, $\theta_l \in \Theta_o$. By Remark \ref{oainclu} and Proposition \ref{homeo}, $F^{(a)}_{\sigma,\nu}(x) -\theta_l=0$ implies $x=H(\theta_l)=\rho_{\theta_l}$. If $l\neq j$, $\rho_{\theta_l}$ does not belong to $[\rho_{\theta_j}-2\delta; \rho_{\theta_j}+2\delta]$ and  the proof of Lemma \ref{pourholomorphe} (1) is complete. 
\end{itemize}

\noindent{Proof of (2):}\\
~~\\
First, note that $0 \notin [\rho_{\theta_j}-2\delta; \rho_{\theta_j}+2\delta].$

\begin{itemize}
\item Assume $\theta_l \notin \Theta_o$.
First,
 if ${\cal Z}^{'} (\frac{1}{\theta_l}) >0$, according to Remark \ref{pasdansimagedeg}, $\frac{1}{\theta_l}$ does not belong to $g_{\tau_{c,\nu}}(\mathbb{R}\setminus\{0\})$ so that the conclusion of Lemma \ref{pourholomorphe} (2) is  true.\\
Now assume that ${\cal Z}^{'} (\frac{1}{\theta_l}) =0$.
By  Remark \ref{compositionlimite}, $1- \theta_l g_{\tau_{c,\nu}}(x)=0$ with $x$ nonnull implies $x={\cal Z}(\frac{1}{\theta_l}).$
In particular, if ${\cal Z}(\frac{1}{\theta_l})=0$, the conclusion of Lemma \ref{pourholomorphe} (2) is true since $0 \notin 
[\rho_{\theta_j}-2\delta; \rho_{\theta_j}+2\delta].$ Hence, in the following, we will deal 
with $\theta_l$ such that  ${\cal Z}(\frac{1}{\theta_l})\neq 0$.
 For any  $u \in {\cal O}^{(m)}$,
we have ${\cal Z}(\frac{1}{\theta_l})
\neq {\cal Z}(\frac{1}{u})$. Indeed, for any $u \in {\cal O}^{(m)}$, there exists $u_1, u_2$ such that $[u_1;u_2] \subset {\cal O}^{(m)}$ and $u_1< u < u_2$. Since  $x \mapsto {\cal Z}(\frac{1}{x})$ is globally nondecreasing on  $\{v\neq 0, \, v \in ^c \mbox{supp~}(\nu),   \, {\cal Z}^{'} (\frac{1}{v}) \leq 0\}$
(see  Remark \ref{croissance globaleS}), we have if $\theta_l <u$, 
 ${\cal Z}(\frac{1}{\theta_l})\leq {\cal Z}(\frac{1}{u_1})< {\cal Z}(\frac{1}{u}),$
and if  $\theta_l >u$, ${\cal Z}(\frac{1}{u})< {\cal Z}(\frac{1}{u_2})\leq {\cal Z}(\frac{1}{\theta_l})$.
It follows, according to  (\ref{car2}), that
${\cal Z}(\frac{1}{\theta_l}) $ belongs to  supp $(\tau_{c,\nu})$ and therefore cannot belong to $[\rho_{\theta_j}-2\delta; \rho_{\theta_j}+2\delta]$ so that the conclusion of Lemma \ref{pourholomorphe} (2)  is true. 

\item Now,
let us consider  $\theta_l \in \Theta_o$. By  Remark \ref{compositionlimite}, $1- \theta_l g_{\tau_{c,\nu}}(x)=0$ with $x$ nonnull implies $x={\cal Z}(\frac{1}{\theta_l})=\rho_{\theta_l}.$ If $l\neq j$, $\rho_{\theta_l}$ does not belong to $[\rho_{\theta_j}-2\delta; \rho_{\theta_j}+2\delta]$ and  the proof of Lemma \ref{pourholomorphe} (2) is complete. $\Box$
\end{itemize}
Let $\theta_j$ be in $\Theta_o$. According to Lemma \ref{pourholomorphe},   $\phi_l$ is an analytic function on $]\rho_{\theta_j}-2\delta; \rho_{\theta_j}+2\delta[\times \mathbb{R}$ for $l\neq j$ and $\phi_j$ is an analytic function on $]\rho_{\theta_j}-2\delta; \rho_{\theta_j}+2\delta[\times \mathbb{R}\setminus \{(\rho_{\theta_j},0)\}$.
Moreover, for any $l$,  $\overline{\phi_l(z)}=\phi_l(\overline{z})$.
We have
\begin{eqnarray*}\frac{1}{\pi}\int \Im \phi_l(t+iy) h_{\delta,j}(t) dt&=& 
\frac{1}{2i\pi}\int_{\rho_{\theta_j} -\delta}^{ \rho_{\theta_j} -\frac{\delta}{2}}h_{\delta,j}(t)\left[\phi_l(t+iy)
-\phi_l(t-iy)\right] dt\\
&&+\frac{1}{2i\pi}\int^{\rho_{\theta_j}+\delta}_{ \rho_{\theta_j} +\frac{\delta}{2}}h_{\delta,j}(t)\left[\phi_l(t+iy)
-\phi_l(t-iy)\right] dt\\
&&+\frac{1}{2i\pi}\int_{\rho_{\theta_j} -\frac{\delta}{2}}^{ \rho_{\theta_j} +\frac{\delta}{2}}\left[\phi_l(t+iy)
-\phi_l(t-iy)\right] dt\\
&=&\Delta_1+ \Delta_2 +\Delta_3.
\end{eqnarray*}
We immediately get  that $\lim_{y \rightarrow 0^+} \Delta_1=0$ and $\lim_{y \rightarrow 0^+} \Delta_2=0$.
Now, \begin{eqnarray*}
\Delta_3&=&\frac{1}{2i\pi}\int_{\rho_{\theta_j} -\frac{\delta}{2}}^{ \rho_{\theta_j} +\frac{\delta}{2}}\left[\phi_l(t+iy)
-\phi_l(t-iy)\right] dt\\
&=& \frac{1}{2i\pi}\int_{\gamma_{j,y,\delta}}\phi_l(z)dz\\
&&-\frac{1}{2\pi}\int_{-y}^{ y}\phi_l(\rho_{\theta_j} -\frac{\delta}{2}+iu)du \\
&&+\frac{1}{2\pi}\int_{-y}^{ y}\phi_l(\rho_{\theta_j} +\frac{\delta}{2}+iu)du,\\&=&\Delta_{3,1}+ \Delta_{3,2} +\Delta_{3,3}
\end{eqnarray*} 
where $\gamma_{j,y,\delta}$ is the clockwise oriented rectangular with corners $\rho_{\theta_j} -\frac{\delta}{2}-iy$,
$\rho_{\theta_j} -\frac{\delta}{2}+iy$, $\rho_{\theta_j} +\frac{\delta}{2}+iy$ and $\rho_{\theta_j} +\frac{\delta}{2}-iy$.
We immediately get   that $\lim_{y \rightarrow 0^+}\Delta_{3,2}=0$ and $\lim_{y \rightarrow 0^+} \Delta_{3,3}=0$.
Moreover, for all  $y$, $$\Delta_{3,1}=\frac{1}{2i\pi}\int_{\gamma_{j,y,\delta}}\phi_l(z)dz =-\mbox{Res~}(\phi_l, \rho_{\theta_j})$$
with $$\mbox{Res~}(\phi_l, \rho_{\theta_j})=0 \mbox{~if~}
l \neq j,$$
\noindent and
$$\mbox{Res~}(\phi_j, \rho_{\theta_j})=\left\{\begin{array}{ll}
\frac{k_j}{F^{(a)'}_{\sigma,\nu}(\rho_{\theta_j})}\, \hspace*{1cm} \mbox{~if~} \, M_N=M_N^W\\
=\frac{k_j \rho_{\theta_j}}{\theta_j  F^{(m)'}_{c,\nu}(\frac{1}{\rho_{\theta_j}})}\, \mbox{~if~} \, M_N=M_N^S.\end{array} \right.$$

\noindent Thus 
\begin{equation}\label{residu}\lim_{y \rightarrow 0^+} \frac{1}{\pi} \int h_{\delta, j}(t) \Im  \phi_l(t+iy)dt 
=0 \mbox{~if~} l\neq j,\end{equation}
and  \begin{equation}\label{residu2}\lim_{y \rightarrow 0^{+}}  \frac{1}{\pi} \int h_{\delta, j}(t)\Im \phi_j(t+iy) dt 
=\left\{\begin{array}{ll}
-k_j H'(\theta_j)\, \hspace*{1cm} \mbox{~if~} \, M_N=M_N^W\\
k_j \frac{{\cal Z}^{'}(\frac{1}{\theta_j})}{\theta_j {\cal Z}(\frac{1}{\theta_j}) } \, \hspace*{0.5cm}\mbox{~if~} \, M_N=M_N^S.\end{array} \right. \end{equation}
Finally from (\ref{ecrit}), (\ref{reste}) (\ref{residu2}) and (\ref{residu}) we deduce that 
$$\lim_{N\rightarrow +\infty} \lim_{y\rightarrow 0^{+}}\frac{1}{\pi} \Im \int  \mathbb{E} \left( \Tr \left[G_{M_N}(t+iy)f_{\eta,l}(A_N)\right] \right)h_{\delta, j}(t) dt=0 \mbox{~if~} l\neq j,$$
and \\

$\lim_{N\rightarrow +\infty} \lim_{y\rightarrow 0^+}\frac{1}{\pi} \Im \int  \mathbb{E} \left( \Tr \left[G_{M_N}(t+iy)f_{\eta,j}(A_N)\right] \right)h_{\delta, j}(t) dt $
$$=\left\{\begin{array}{ll}
-k_j H'(\theta_j)\, \hspace*{1cm} \mbox{~if~} \, M_N=M_N^W\\
k_j \frac{{\cal Z}^{'}(\frac{1}{\theta_j})}{\theta_j {\cal Z}(\frac{1}{\theta_j}) } \, \hspace*{0.5cm} \mbox{~if~} \, M_N=M_N^S.\end{array} \right. $$
Then, Proposition \ref{cvev} follows by Proposition \ref{premetap} and Lemma \ref{approxcauchy}.$\Box$

\section{Technical results specific to each model}
\begin{lemme}\label{Lipschitz}
~~\\
\begin{itemize}
 \item[(1)] For any $N \times N$ Hermitian matrix $X$,\\
~
$\left\{\left(X(i,i)\right)_{1\leq i\leq N}\,\left(\sqrt{2}\Re X(i,j)\right)_{1\leq i<j \leq N} \,\left(\sqrt{2}\Im X(i,j)\right)_{1\leq i< j \leq N}\right\}$ $$ \mapsto \Tr \left[ h_{\delta,j} (X+A_N) f_{\eta,l}(A_N)\right]$$
\noindent  is Lipschitz with constant bounded by $\sqrt{k_l} \Vert h_{\delta,j} \Vert_{Lip}$.
\item[(2)]For any $N \times p$  matrix $B$,
$$\left\{\left(\Re B(i,j),~ \Im B(i,j)\right)_{1\leq i \leq N, 1 \leq j \leq p}\right\} \mapsto \Tr \left[ h_{\delta,j} (A_N^{\frac{1}{2}}BB^*A_N^{\frac{1}{2}}) f_{\eta,l}(A_N)\right]$$
\noindent  is Lipschitz with constant bounded by $\sqrt{k_l} \sqrt{2C} \Vert \tilde{h}_{\delta,j} \Vert_{Lip}$
where $\tilde{h}_{\delta,j} (x) =h_{\delta,j}(x^2)$ and $C=\sup_N \Vert A_N \Vert $.
\end{itemize}
\end{lemme}
\noindent {\bf Proof} Given two $N \times N$ Hermitian matrices $X$ and $X^{'}$,
we have (using Lemma \ref{extlipschitz}) that \\

$\left|\Tr \left[ h_{\delta,j} (X+A_N) f_{\eta,l}(A_N)\right]-\Tr \left[ h_{\delta,j} (X^{'}+A_N) f_{\eta,l}(A_N)\right]\right|$
$$\leq \Vert f_{\eta,l}(A_N)\Vert_2 \Vert X-X^{'} \Vert_2 \Vert h_{\delta,j} \Vert_{Lip}$$
and (1) follows since $\Vert f_{\eta,l}(A_N)\Vert_2 =\sqrt{k_l}$.\\
To prove (2) we will make use of a useful observation already made in \cite{GuZei}.
Let us introduce the $(N+p) \times (N+p)$ matrices $${\cal M}_{N+p}(B) =\left( \begin{array}{ll} 0_{p\times p}~~~ B^*A_N^{\frac{1}{2}}\\ A_N^{\frac{1}{2}}B ~~~ 0_{N\times N} \end{array} \right),$$    
$${\cal N}_{N+p}=\left( \begin{array}{ll} 0_{p\times p}~~~ 0_{p\times N}\\ 0_{N\times p} ~~~ f_{\eta,l}(A_N)\end{array} \right).$$
It is easy to see that 
$$\Tr \left[ h_{\delta,j} (A_N^{\frac{1}{2}}BB^*A_N^{\frac{1}{2}}) f_{\eta,l}(A_N)\right]=
  \Tr \left[ \tilde{h}_{\delta,j} ({\cal M}_{N+p}(B)) {\cal N}_{N+p}\right],$$
  where $\tilde{h}_{\delta,j}(x) =h_{\delta,j}(x^2)$. Note that since $h_{\delta,j}$ is a ${\cal C}^\infty$ compactly supported function,
  $\tilde{h}_{\delta,j}$ is obviously a Lipschitz function.
  Hence \\
  
  $\left| \Tr \left[ h_{\delta,j} (A_N^{\frac{1}{2}}BB^*A_N^{\frac{1}{2}}) f_{\eta,l}(A_N)\right]
  -\Tr \left[ h_{\delta,j} (A_N^{\frac{1}{2}}B^{'}(B^{'})^*A_N^{\frac{1}{2}}) f_{\eta,l}(A_N)\right]\right|$
 \begin{eqnarray}&= &\left| \Tr \left[ \tilde{h}_{\delta,j} ({\cal M}_{N+p}(B)) {\cal N}_{N+p}\right]- \Tr \left[ \tilde{h}_{\delta,j} ({\cal M}_{N+p}(B^{'})) {\cal N}_{N+p}\right]\right| \nonumber\\ &\leq &
  \left\| {\cal N}_{N+p}\right\|_2 \left\| \tilde{h}_{\delta,j} ({\cal M}_{N+p}(B))-
 \tilde{h}_{\delta,j} ({\cal M}_{N+p}(B^{'}))\right\|_2 \nonumber\\
 &\leq & \sqrt{k_l}\left\|\tilde{h}_{\delta,j} \right\|_{Lip} \left\| {\cal M}_{N+p}(B)-{\cal M}_{N+p}(B^{'}) \right\|_2.\label{lip1}
 \end{eqnarray}
 \noindent where we used Lemma \ref{extlipschitz} in the last line.
 Now, \begin{eqnarray}\left\| {\cal M}_{N+p}(B)-{\cal M}_{N+p}(B^{'}) \right\|^2_2&=& 2 \Tr \left[ (B-B^{'})^{*}A_N (B-B^{'}) \right]\nonumber \\ &\leq & 2\left\| A_N \right\| \left\| B-B^{'}\right\|^2_2
\nonumber \\ & \leq & 2 C  \left\| B-B^{'}\right\|^2_2.\label{lip2}\end{eqnarray}
 
\noindent (2) readily follows from (\ref{lip1}) and (\ref{lip2}). 
$\Box$\\

We have proved in Lemma 3.3 \cite{CDFF10} that $\forall z \in \mathbb{C}\setminus \mathbb{R}$,
\begin{equation}\label{estimfondakW}
\E (\tilde G^W_{kk}(z)) = \frac{1}{(z-\sigma ^2g^W_N(z)-\gamma_k)}
+R_{k,N}(z),
\end{equation}
with $$\vert R_{k,N}(z) \vert \leq \frac{P(\vert \Im z \vert^{-1})}{N}$$
\noindent for some
polynomial $P$ with nonnegative coefficients.
Note that \begin{equation} \label{imageW} \vert \Im (z-\sigma ^2g^W_N(z)-\gamma_k)\vert \geq \vert \Im z \vert.\end{equation}
We are going to establish the following  similar result for the sample covariance matrix setting using many ideas from \cite{BaiSil06}.

\begin{proposition} \label{estimresolvwish} There exists a polynomial $P$ with nonnegative coefficients and a sequence $a_N$ of nonnegative
real numbers converging to zero when $N$ goes to infinity such that, for any $k$ in $\{1, \ldots,N\}$, any $z$ in  $\mathbb{C}\setminus \mathbb{R}$,
\begin{equation} \label{submat}
\E (\tilde G^S_{kk}(z)) = \frac{1}{z-\gamma_k(1-\frac{N}{p}+\frac{N}{p}zg^S_N(z))}
+R_{k,N}(z),
\end{equation}
with $$\left|R_{k,N}(z) \right| \leq (1+\vert z\vert)^2 P(\vert \Im z \vert^{-1})a_N.$$
\end{proposition}

\noindent {\bf Proof}
Let $J$ be a $N\times N$ matrix and $u$ be a vector in $\mathbb{C}^N$ such that $J$ and $J +uu^*$ are invertible  then
\begin{eqnarray*}
u^* J^{-1} u u^* \left( J+ uu^* \right)^{-1} & = & u^* J^{-1} \left(u u^* +J \right) \left( J+ uu^* \right)^{-1} 
- u^* \left( J+ uu^* \right)^{-1}\\&=& u^* J^{-1}-u^* \left( J+ uu^* \right)^{-1}
\end{eqnarray*}
so that \begin{equation}\label{trick}
u^* \left( J+ uu^* \right)^{-1}=\frac{u^* J^{-1}}{1+u^* J^{-1}u}.
\end{equation}
Hence if $u_1,\ldots,u_p$ are $p$ vectors and $X=\sum_{i=1}^p u_i u_i^*$, denoting by $G_X(z)$ the resolvent
$(zI-X)^{-1}$,
(\ref{trick}) yields that for any $i\in \{1,\ldots,p\}$, for any $z\in \mathbb{C}\setminus \mathbb{R}$,
\begin{equation}\label {etape1}u_iG_X(z) = \frac{u_i^* \left(zI-\sum_{l\neq i} u_l u_l^*\right)^{-1}}{1-u_i^* \left(zI-\sum_{l\neq i} u_l u_l^*\right)^{-1}u_i}.\end{equation}
Multiplying (\ref{etape1}) by $u_i$ and summing in $i$ yields
\begin{equation}\label{etape2}
XG_X(z) =\sum_{i=1}^p \frac{u_i u_i^* \left(zI-\sum_{l\neq i} u_l u_l^*\right)^{-1}}{1-u_i^* \left(zI-\sum_{l\neq i} u_l u_l^*\right)^{-1}u_i}.
\end{equation}
From (\ref{etape2}) and the resolvent identity
$$-I+zG_X(z) =XG_X(z),$$ we deduce that 
\begin{equation}\label{Gkk}
\left(G_X(z)\right)_{kk} =\frac{1}{z} +
\sum_{i=1}^p \frac{\left[u_i u_i^* \left(zI-\sum_{l\neq i} u_l u_l^*\right)^{-1}\right]_{kk}}{z \left[1-u_i^* \left(zI-\sum_{l\neq i} u_l u_l^*\right)^{-1}u_i\right]}.
\end{equation}
Noticing that \begin{equation}\label{remarque} \tilde{G}^S(z)= G_{\frac{1}{p} D^{\frac{1}{2}}U X_N^S U^* D^{\frac{1}{2}}}\end{equation}
and that $$\frac{1}{p} D^{\frac{1}{2}}U X_N^S U^* D^{\frac{1}{2}} =\sum_{i=1}^p u_i u_i^*$$
where $u_i =\frac{1}{\sqrt{p}} D^{\frac{1}{2}} U x_i$ and $x_i$ is the $i$th column of $B_N$,
we deduce from (\ref{Gkk}) that 
\begin{equation}\label{Gkk2}
\tilde{G}^S_{kk}(z)=\frac{1}{z} + \frac{1}{p} \sum_{i=1}^p 
\frac{\left[D^{\frac{1}{2}} U x_i x_i^* U^* D^{\frac{1}{2}}  \left(zI-\frac{1}{p}\sum_{l\neq i} D^{\frac{1}{2}} U x_l x_l^* U^* D^{\frac{1}{2}}\right)^{-1}\right]_{kk}}{z \left\{1-\frac{1}{p} x_i^* U^* D^{\frac{1}{2}}  \left(zI- \frac{1}{p} \sum_{l\neq i} D^{\frac{1}{2}} U x_l x_l^* U^* D^{\frac{1}{2}}\right)^{-1}D^{\frac{1}{2}} U x_i \right\}}.
\end{equation}
Set for $i=1,\ldots,p$, $$y_i = \frac{1}{\sqrt{p}} A_N^{\frac{1}{2}}x_i$$
and $$M_N^{(i)} = \frac{1}{{p}} \sum_{l\neq i} A_N^{\frac{1}{2}}x_l x_l^* A_N^{\frac{1}{2}}=\sum_{l\neq i} y_l y_l^*.$$
Note that the $y_i$'s are i.i.d and that $y_i$ is independent of $M_N^{(i)}$.
Note also that $$M_N^S =M_N^{(i)} + y_i y_i^*.$$
(\ref{Gkk2}) can be rewritten as follows
\begin{equation}\label{Gkk3}
\tilde{G}^S_{kk}(z)=\frac{1}{z} + \frac{1}{p} \sum_{i=1}^p 
\frac{\left[D^{\frac{1}{2}} U x_i x_i^*  A_N^{\frac{1}{2}}  \left(zI-M_N^{(i)}\right)^{-1}U^*\right]_{kk}}{z \left\{1-y_i^*   \left(zI- M_N^{(i)}\right)^{-1} y_i \right\}}
\end{equation}
Now applying (\ref{Gkk}) with $u_i= y_i$ and $X= M_N^S$, summing on $k$ and dividing by $N$ we have
\begin{eqnarray}
\tr_N\left( G_{M_N^S} (z)\right) &=& \frac{1}{z} + \frac{1}{N} \sum_{i=1}^p 
\frac{y_i^* \left(zI-M_N^{(i)}\right)^{-1}y_i}{z \left\{1-y_i^*   \left(zI- M_N^{(i)}\right)^{-1} y_i \right\}}\nonumber\\
&=& 
\frac{1}{z}-\frac{p}{N} \frac{1}{z} + \frac{1}{N} \sum_{i=1}^p 
\frac{1}{z \left\{1-y_i^*   \left(zI- M_N^{(i)}\right)^{-1} y_i \right\}}.\label{trace}
\end{eqnarray}
Let us define the $p\times p$ matrix 
$$\underline{M}_N^S= \frac{1}{p}B_N^* A_N B_N.$$
Since 
\begin{equation}\label{baroupas}
\tr_p\left( G_{\underline{M}_N^S} (z)\right)=\frac{1-\frac{N}{p}}{z} +\frac{N}{p}
\tr_N\left( G_{M_N^S} (z)\right),\end{equation}
we deduce from (\ref{trace}) that 
\begin{equation}\label{tracebar}
\tr_N\left( G_{\underline{M}_N^S} (z)\right) = \frac{1}{p} \sum_{i=1}^p 
\frac{1}{z \left\{1-y_i^*   \left(zI- M_N^{(i)}\right)^{-1} y_i \right\}}.
\end{equation}
Following the ideas of Section 6.4.1 of \cite{BaiSil06}, we are going to establish the following preliminary lemma.
\begin{lemme} \label{normeL2}
There exists a constant $K>0$ and a sequence of nonnegative numbers $a_N$ converging to zero when $N$ goes to infinity such that for each $i=1, \ldots,p$, $\forall z \in \mathbb{C}\setminus \mathbb{R}$,
$$\left\| \tr_p\left( G_{\underline{M}_N^S} (z)\right) - \frac{1}{z \left\{1-y_i^*   \left(zI- M_N^{(i)}\right)^{-1} y_i \right\}} \right\|_{L^2} \leq K \frac{\vert z \vert }{\vert \Im z \vert^3}a_N.$$
\end{lemme}
\noindent We have from (\ref{tracebar})\\

$\tr_p\left( G_{\underline{M}_N^S} (z)\right) - \frac{1}{z \left\{1-y_i^*   \left(zI- M_N^{(i)}\right)^{-1} y_i \right\}} $
$$= \frac{1}{p} \sum_{l\neq i}
\frac{y_l^*   \left(zI- M_N^{(l)}\right)^{-1} y_l -  y_i^*   \left(zI- M_N^{(i)}\right)^{-1} y_i   }{z \left\{1-y_l^*   \left(zI- M_N^{(l)}\right)^{-1} y_l \right\}\left\{1-y_i^*   \left(zI- M_N^{(i)}\right)^{-1} y_i \right\}}.$$

\begin{lemme}\label{maj}
For   any $N\times N$ positive semidefinite matrix $H$, any vector $v$ in $\mathbb{C}^N$ and any $z$ in $
\mathbb{C}\setminus \mathbb{R}$,
$$\frac{1}{\vert z \left\{1-v^*   (zI- H)^{-1} v \right\}\vert }  \leq \frac{1}{\vert \Im z \vert }.$$
\end{lemme}
\noindent {\bf Proof:}\begin{eqnarray*}
 \Im \left\{ z v^*  (zI- H)^{-1} v \right\}&=& \frac{1}{2i} \left\{
 z v^*  (zI- H)^{-1} v -\bar{z} v^* (\bar{z}I- H)^{-1} v \right\}\\
 &=&
  \frac{1}{2i} v^*\left\{
 (I- \frac{1}{z}H)^{-1}  -(I- \frac{1}{\bar{z}} H)^{-1} \right\}v\\
&=&
  -\frac{\Im z}{\vert z \vert^2} v^*
 (I- \frac{1}{z}H)^{-1}H(I- \frac{1}{\bar{z}} H)^{-1} v.
 \end{eqnarray*}
 Hence \begin{eqnarray*}
 \vert \Im \left[ z \left\{1-v^*   (zI- H)^{-1} v \right\} \right] \vert &=&
 \vert \Im z \vert \left\{ 1 + \frac{\1}{\vert z \vert^2}v^*
 (I- \frac{1}{z}H)^{-1}H(I- \frac{1}{\bar{z}} H)^{-1} v \right\}\\
 & \geq & \vert \Im z \vert
 \end{eqnarray*}
 and Lemma \ref{maj} follows.$\Box$\\
 
 \noindent According to Lemma \ref{maj}, for any $l$ and $i$ in $\{1, \ldots,p\}$,
 $$\frac{1}{\left|z \left\{1-y_l^*   \left(zI- M_N^{(l)}\right)^{-1} y_l \right\}\left\{1-y_i^*   \left(zI- M_N^{(i)}\right)^{-1} y_i \right\}\right| }\leq \frac{\vert z \vert}{\vert \Im z \vert^2 }.$$
 Hence \\
 
 $\left\| \tr_p\left( G_{\underline{M}_N^S} (z)\right) - \frac{1}{z \left\{1-y_i^*   \left(zI- M_N^{(i)}\right)^{-1} y_i \right\}} \right\|_{L^2}$
 $$\leq \frac{\vert z \vert}{\vert \Im z \vert^2 } \frac{1}{p} \sum_{l\neq i} \left\| 
y_l^*   \left(zI- M_N^{(l)}\right)^{-1} y_l -  y_i^*   \left(zI- M_N^{(i)}\right)^{-1} y_i \right\|_{L^2}.$$
We have\\
 
$
y_l^*   \left(zI- M_N^{(l)}\right)^{-1} y_l -  y_i^*   \left(zI- M_N^{(i)}\right)^{-1} y_i$
\begin{eqnarray*}&=&
y_l^*   \left(zI- M_N^{(l)}\right)^{-1} y_l -  \frac{1}{p}
\Tr \left[\left(zI- M_N^{(l)}\right)^{-1} A_N \right]\\ & &
+\frac{1}{p}
\Tr \left[\left(zI- M_N^{(l)}\right)^{-1} A_N \right] - \frac{1}{p}
\Tr \left[\left(zI- M_N^{(i)}\right)^{-1} A_N \right]\\
&& +\frac{1}{p}
\Tr \left[\left(zI- M_N^{(i)}\right)^{-1} A_N \right]
-y_i^*   \left(zI- M_N^{(i)}\right)^{-1} y_i\\
&=& \Delta_l + \Delta_{l,i} - \Delta_i.
\end{eqnarray*}
We have 
\begin{eqnarray*}
\vert \Delta_{l,i} \vert &\leq & \left| \frac{1}{p}
\Tr \left[\left\{\left(zI- M_N^{(l)}\right)^{-1}-\left(zI- M_N^S \right)^{-1} \right\}A_N \right] \right| \\
& & +\left| \frac{1}{p}
\Tr \left[\left\{\left(zI- M_N^S\right)^{-1}-\left(zI- M_N^{(i)} \right)^{-1} \right\}A_N \right] \right|\\
&\leq & \frac{2 \Vert A_N \Vert}{\vert \Im z \vert p}\leq \frac{2 C}{\vert \Im z \vert p} 
\end{eqnarray*}
where we used Lemma 6.9 \cite{BaiSil06}) in the last line. \\
Now we have for any $l=1,\ldots,p$,
$$y_l^*   \left(zI- M_N^{(l)}\right)^{-1} y_l= \frac{1}{p} x_l^* A_N^{\frac{1}{2}}\left(zI- M_N^{(l)}\right)^{-1} 
A_N^{\frac{1}{2}}x_l$$
so that according to Proposition \ref{QF} in the Appendix,
for any $l=1,\ldots,p$, $$\Vert \Delta_{l} \Vert_{L^2} \leq \sqrt{K} \sqrt{\frac{N}{p}} \frac{1}{\sqrt{p}} 
 \frac{ \Vert A_N \Vert}{\vert \Im z \vert }
 \leq \sqrt{\frac{N}{p}} \frac{1}{\sqrt{p}} 
 \frac{ C}{\vert \Im z \vert }.$$
 It follows that 
 $$\left\| \tr_p\left( G_{\underline{M}_N^S} (z)\right) - \frac{1}{z \left\{1-y_i^*   \left(zI- M_N^{(i)}\right)^{-1} y_i \right\}} \right\|_{L^2} \leq \frac{2C \vert z \vert }{\vert \Im z \vert^3}
 \left( \frac{1}{p}+ \sqrt{\frac{N}{p}} \frac{1}{\sqrt{p}} \right)$$
 and the proof of Lemma \ref{normeL2} is complete. $\Box$
 \begin{lemme}\label{phi}
 There exists a constant $K>0$ such that, for any $k=1, \ldots, N$, for any $i=1,\ldots,p$ and any $z$ in $\mathbb{C}\setminus \mathbb{R}$,
 $$\left\|
 \left[D^{\frac{1}{2}} U x_i x_i^*  A_N^{\frac{1}{2}}  \left(zI-M_N^{(i)}\right)^{-1}U^*\right]_{kk}
 \right\|_{L^2} \leq \frac{K}{\vert \Im z \vert}.$$
 \end{lemme}
 \noindent {\bf Proof:}
 Note that
 \begin{eqnarray*}
 \left[D^{\frac{1}{2}} U x_i x_i^*  A_N^{\frac{1}{2}}  \left(zI-M_N^{(i)}\right)^{-1}U^*\right]_{kk}&=&
 \Tr D^{\frac{1}{2}} U x_i x_i^*  A_N^{\frac{1}{2}}  \left(zI-M_N^{(i)}\right)^{-1}U^*E_{kk}\\
 &=&  x_i^*  A_N^{\frac{1}{2}}  \left(zI-M_N^{(i)}\right)^{-1}U^*E_{kk}D^{\frac{1}{2}} U x_i.
 \end{eqnarray*}
 Thus, according to Proposition \ref{QF} in the Appendix,\\
 
 $
\left\|
 \left[D^{\frac{1}{2}} U x_i x_i^*  A_N^{\frac{1}{2}}  \left(zI-M_N^{(i)}\right)^{-1}U^*\right]_{kk}
 -\Tr A_N^{\frac{1}{2}}  \left(zI-M_N^{(i)}\right)^{-1}U^*E_{kk}D^{\frac{1}{2}} U 
 \right\|_{L^2}$ $$\leq K
 \left[ \Tr A_N^{\frac{1}{2}}  \left(zI-M_N^{(i)}\right)^{-1}U^*E_{kk}D E_{kk} U  \left(\bar{z}I-M_N^{(i)}\right)^{-1} A_N^{\frac{1}{2}}\right]^{\frac{1}{2}} \leq \frac{C}{\vert \Im z\vert}.
 $$
 Since moreover 
 $$\left| \Tr A_N^{\frac{1}{2}}  \left(zI-M_N^{(i)}\right)^{-1}U^*E_{kk}D^{\frac{1}{2}} U \right| \leq \frac{C}{\vert \Im z\vert},$$
 we deduce that $$\left\|
 \left[D^{\frac{1}{2}} U x_i x_i^*  A_N^{\frac{1}{2}}  \left(zI-M_N^{(i)}\right)^{-1}U^*\right]_{kk}
\right\|_{L^2} \leq \frac{2C}{\vert \Im z \vert}.$$
 \hspace*{10cm}$\Box$\\
 We will need this last lemma concerning the variance of 
 $\tr_p\left( G_{\underline{M}_N^S} (z)\right)$.
 \begin{lemme}\label{variancebar} There exists some polynomial $P$ with nonnegative coefficients  such that, $\forall z \in \mathbb{C}\setminus \mathbb{R}$,
 $$\left\| \tr_p\left( G_{\underline{M}_N^S} (z)\right) - \mathbb{E}\left( \tr_p\left( G_{\underline{M}_N^S} (z)\right)\right) \right\|_{L^2} \leq \frac{1}{p^2} \left( \vert z\vert +1\right)^2 P\left(\vert \Im z \vert^{-1}\right).$$
 \end{lemme}
 \noindent{\bf Proof:}
 Let us define $\Psi: \R ^{2(p\times N)} \rightarrow {\mathcal M }_{N\times p}(\C)$ by 
$$\Psi:~~\{x_{ij},y_{ij},i=1,\ldots,N,j=1,\ldots,p\}\rightarrow 
\sum_{i=1,\ldots,N}\sum_{j=1,\ldots,p} \left( x_{ij} +\sqrt{-1} y_{ij} \right) E_{ij}.$$
Let $F$ be a smooth complex function on ${ M}_{N\times p}(\C)$ and 
define the complex function $f$ on  $\R^{2(p\times N)}$ by setting $f=F\circ \Psi$.
Then,$$\Vert {\rm grad} f(u)\Vert = \sup_{V\in { M }_{N\times p}(\C), Tr VV^*=1}
\left| \frac{d}{dt} F(\Psi(u)+tV)\vert _{t=0}\right| .$$
We have $B_N=\Psi(\Re ((B_N)_{ij}),\Im ((B_N)_{ij}),1\leq i\leq N,1\leq j\leq p)$ 
where the distribution of $\{\sqrt{2}\Re ((B_N)_{ij}),\sqrt{2}\Im ((B_N)_{ij}),1\leq i\leq N,1\leq j\leq p\}$ 
satisfies a Poincar\'e inequality with constant $C_{PI}$.\\
Hence consider $F:~B \rightarrow \tr _N(zI_N -A_N^{\frac{1}{2}} \frac{B B^*}{p}A_N^{\frac{1}{2}} )^{-1}$.\\
Let $ V\in { M }_{N\times p}(\C)$ such that $ Tr VV^*=1$.
\begin{eqnarray}
\frac{d}{dt} F(B+tV)\vert _{t=0}&=&\frac{1}{Np} \Tr(G_{M_N^S}(z) A_N^{\frac{1}{2}} V B^* A_N^{\frac{1}{2}} G_{M_N^S}(z)) \nonumber \\
&&+ \frac{1}{Np} \Tr(G_{M_N^S}(z) A_N^{\frac{1}{2}} B V^* A_N^{\frac{1}{2}} G_{M_N^S}(z)). \label{derive}
\end{eqnarray}
By Cauchy-Schwartz inequality, we have \\

$
\left| \frac{1}{Np} \Tr(G_{M_N^S}(z) A_N^{\frac{1}{2}} V B_N^* A_N^{\frac{1}{2}} G_{M_N^S}(z)) \right| $
\begin{eqnarray*}
&\leq &\frac{1}{N \sqrt{p}} 
\left[\Tr \frac{B_N^*A_N^{\frac{1}{2}} \left[G_{M_N^S}(z)\right]^2 A_N \left[G_{M_N^S}(\bar{z})\right]^2 A_N^{\frac{1}{2}} B_N }{p}\right] ^{\frac{1}{2}}
(\Tr VV^* )^{\frac{1}{2}}\\
&= & \frac{1}{\sqrt{Np}}\left[\tr_N M_N^S \left[G_{M_N^S}(z)\right]^2 A_N \left[G_{M_N^S}(\bar{z})\right]^2 
\right] ^{\frac{1}{2}}.
\end{eqnarray*}
Since by the resolvent identity
$$M_N^S G_{M_N^S}(z)=-I_N +z G_{M_N^S},$$ we have \\

$\tr_N M_N^S \left[G_{M_N^S}(z)\right]^2 A_N \left[G_{M_N^S}(\bar{z})\right]^2$ $$ =
-\tr_N G_{M_N^S}(z) A_N \left[G_{M_N^S}(\bar{z})\right]^2 +z \tr_N  \left[G_{M_N^S}(z)\right]^2 A_N \left[G_{M_N^S}(\bar{z})\right]^2,$$
we can deduce  that  
$$\tr_N M_N^S \left[G_{M_N^S}(z)\right]^2 A_N \left[G_{M_N^S}(\bar{z})\right]^2 \leq \frac{\Vert A_N \Vert}{\vert \Im z \vert^3} +\frac{\vert z \vert \Vert A_N \Vert}{\vert \Im z \vert^4}\leq C(\vert z \vert +1) P\left( \vert \Im z \vert^{-1} \right),$$
\noindent where $P$ is a polynomial with non negative coefficients.
 Hence 
 $$\left| \frac{1}{Np} \Tr(G_{M_N^S}(z) A_N^{\frac{1}{2}} V B_N^* A_N^{\frac{1}{2}} G_{M_N^S}(z)) \right| \leq
 \frac{1}{\sqrt{Np}}C(\vert z \vert +1) P\left( \vert \Im z \vert^{-1} \right).$$
Since a similar upper bound can be obtained in the same way for the second term on the right hand side of (\ref{derive})
we deduce that 
$$\mathbb{E}\left( \left( \sup_{V\in { M }_{N\times p}(\C), Tr VV^*=1}
\left| \frac{d}{dt} F(B_N+tV)\vert _{t=0}\right|\right)^2 \right)^{\frac{1}{2}} \leq  \frac{2C}{\sqrt{Np}}(\vert z \vert +1) P\left( \vert \Im z \vert^{-1} \right).$$
 Therefore, Poincar\'e inequality yields 
 $$\left\| \tr_N\left( G_{{M}_N^S} (z)\right) - \mathbb{E}\left( \tr_N\left( G_{{M}_N^S} (z)\right)\right) \right\|_{L^2}\leq
 \frac{1}{{Np}}(\vert z \vert +1)^2 Q\left(\vert  \Im z \vert^{-1} \right)$$
 \noindent where $Q$ is a polynomial with non negative coefficients.
 Now, since 
 $$\tr_p\left( G_{\underline{M}_N^S} (z)\right) - \mathbb{E}\left( \tr_p\left( G_{\underline{M}_N^S} (z)\right)\right) 
=\frac{N}{p} \left[\tr_N\left( G_{{M}_N^S} (z)\right) - \mathbb{E}\left( \tr_N\left( G_{{M}_N^S} (z)\right)\right) \right],$$
Lemma \ref{variancebar} follows. $\Box$\\

\noindent Using (\ref{Gkk3}), we have for any $z$ in $\mathbb{C}\setminus \mathbb{R}$,
\begin{eqnarray*}
\mathbb{E}\left( \tilde{G}^S_{kk}(z)\right) &=&\frac{1}{z} + 
\frac{1}{p} \sum_{i=1}^p \mathbb{E} \left(\Phi_{i,k} \right)  \mathbb{E} \left(\tr_p G_{\underline{M}_N^S} (z)\right)  \\&& +\frac{1}{p} \sum_{i=1}^p \mathbb{E} \left[\Phi_{i,k} \left\{\tr_p G_{\underline{M}_N^S} (z)
- \mathbb{E} \left(\tr_p G_{\underline{M}_N^S} (z)\right)\right\}\right]\\
&& +\frac{1}{p} \sum_{i=1}^p \mathbb{E} \left[\Phi_{i,k} \left\{  \frac{1}{z \left\{1-y_i^*   \left(zI- M_N^{(i)}\right)^{-1} y_i \right\}}   
 -\tr_p G_{\underline{M}_N^S} (z)
\right\}\right]
\end{eqnarray*}
where $$\Phi_{i,k}=
\left[D^{\frac{1}{2}} U x_i x_i^*  A_N^{\frac{1}{2}}  \left(zI-M_N^{(i)}\right)^{-1}U^*\right]_{kk}.$$
By Cauchy-Schwartz inequality, using Lemmas \ref{phi} and \ref{normeL2}, we easily have that there exists a constant $K$
and a sequence of nonnegative numbers $a_N$ converging towards zero when $N$ goes to infinity such that for any $k=1,\ldots,N$,
$$\left| \frac{1}{p} \sum_{i=1}^p \mathbb{E} \left[\Phi_{i,k} \left\{  \frac{1}{z \left\{1-y_i^*   \left(zI- M_N^{(i)}\right)^{-1} y_i \right\}}   
 -\tr_p G_{\underline{M}_N^S} (z)
\right\}\right]\right| \leq \frac{K \vert z \vert}{\vert \Im z \vert^4} a_N.$$
By Cauchy-Schwartz inequality, using Lemmas \ref{phi} and \ref{variancebar}, we also have that there exists a polynomial $P$ with
nonnegative coefficients such that 
$$\left| \frac{1}{p} \sum_{i=1}^p \mathbb{E} \left[\Phi_{i,k} \left\{\tr_p G_{\underline{M}_N^S} (z)
- \mathbb{E} \left(\tr_p G_{\underline{M}_N^S} (z)\right)\right\}\right]\right|\leq \frac{1}{p^2}(\vert z \vert +1)^2
P(\vert \Im z \vert^{-1}).$$
Thus 
\begin{equation}\label{eqinter} \mathbb{E}\left( \tilde{G}^S_{kk}(z)\right) =\frac{1}{z} + 
\frac{1}{p} \sum_{i=1}^p \mathbb{E} \left(\Phi_{i,k} \right)  \mathbb{E} \left(\tr_p G_{\underline{M}_N^S} (z)\right) + \Delta_N(k)\end{equation}
where there exists a polynomial $Q$ with
nonnegative coefficients and a sequence of nonnegative numbers $b_N$ converging towards zero when $N$ goes to infinity
such that, for any $k=1,\ldots,N$,
$$\vert \Delta_N(k) \vert \leq (\vert z \vert +1)^2
Q(\vert \Im z \vert^{-1}) b_N.$$

\noindent Now, one can easily see that 
\begin{eqnarray}
\mathbb{E} \left(\Phi_{i,k} \right) &=&  \gamma_k \mathbb{E} ( [U
(zI-M_N^{(i)})^{-1}U^*]_{kk} ) \nonumber \\
&=&  \gamma_k \mathbb{E} ( [
(zI-\sum_{l\neq i} u_l u_l^*)^{-1}]_{kk} ) \label{eqphi}
\end{eqnarray} 
where $$u_i=\frac{1}{\sqrt{p}}D^{\frac{1}{2}}Ux_i.$$
\begin{lemme}\label{difftilde} 
There exists a polynomial $P$ with
nonnegative coefficients such that for any $i=1, \ldots,p$, any $k=1,\ldots,N$, any
$ z \in \mathbb{C}\setminus \mathbb{R}$,
$$\left|\mathbb{E}\left( \tilde{G}^S_{kk}(z)\right)- \mathbb{E} ( [
(zI-\sum_{l\neq i} u_l u_l^*)^{-1}]_{kk}) \right| \leq \frac{1}{p}(\vert z \vert +1)
P(\vert \Im z \vert^{-1}).$$
\end{lemme}
\noindent {\bf Proof:} Remember that according to (\ref{remarque}),
$$\mathbb{E}( \tilde{G}^S_{kk}(z))=\mathbb{E} ( [
(zI-\sum_{l=1}^p u_l u_l^*)^{-1}]_{kk} ).
$$
By the formula (3.3.4) in \cite{BaiSil06},
we have $$[
(zI-\sum_{l=1}^p u_l u_l^*)^{-1}]_{kk}=[
(zI-\sum_{l\neq i} u_l u_l^*)^{-1}]_{kk}+ \frac{\psi_{i,k} }{1-u_i^*\left(zI-\sum_{l\neq i} u_l u_l^*\right)^{-1}u_i}$$
where $$\psi_{i,k}= \left[
\left(zI-\sum_{l\neq i} u_l u_l^*\right)^{-1}u_i u_i^*\left(zI-\sum_{l\neq i} u_l u_l^*\right)^{-1}\right]_{kk} .$$

\noindent Noticing that $$\psi_{i,k}= \frac{1}{p} x_i^* U^* D^{\frac{1}{2}} (zI-\sum_{l\neq i} u_l u_l^*)^{-1}
E_{kk} (zI-\sum_{l\neq i} u_l u_l^*)^{-1} D^{\frac{1}{2}} U x_i,$$
we have by Proposition \ref{QF} that\\

$
\left\| \psi_{i,k}- \frac{1}{p} \Tr D \left(zI-\sum_{l\neq i} u_l u_l^*\right)^{-1}
E_{kk} \left(zI-\sum_{l\neq i} u_l u_l^*\right)^{-1} \right\|_{L^2}$
 $$\leq  \frac{K}{p}
\left\{ \Tr D (zI-\sum_{l\neq i} u_l u_l^*)^{-1}
E_{kk} (zI-\sum_{l\neq i} u_l u_l^*)^{-1} 
(\bar{z}I-\sum_{l\neq i} u_l u_l^*)^{-1}
E_{kk} (\bar{z}I-\sum_{l\neq i} u_l u_l^*)^{-1} D \right\}^{\frac{1}{2}}$$
$\leq\frac{CK}{p \vert \Im z \vert^{2}}.
$\\
~

\noindent Using also Lemma \ref{maj}, we readily have that \\

$\left|\mathbb{E}\left( \tilde{G}^S_{kk}(z)\right)- \mathbb{E} \left( \left[
\left(zI-\sum_{l\neq i} u_l u_l^*\right)^{-1}\right]_{kk} \right) \right|$
$$\leq \frac{\vert z \vert}{p \vert \Im z \vert} \left\{ \frac{CK}{ \vert \Im z \vert^{2}} +
\mathbb{E} \left( \left| \Tr D(zI-\sum_{l\neq i} u_l u_l^*)^{-1} E_{kk} (zI-\sum_{l\neq i} u_l u_l^*)^{-1} \right| \right) \right\}.$$
Since $$\left| \Tr D(zI-\sum_{l\neq i} u_l u_l^*)^{-1} E_{kk} (zI-\sum_{l\neq i} u_l u_l^*)^{-1} \right| \leq \frac{C}{ \vert \Im z \vert^{2}},$$
Lemma \ref{difftilde} readily follows. $\Box$\\

\noindent Hence (\ref{eqphi}) and Lemma \ref{difftilde} yield
$$\mathbb{E}\left( \Phi_{i,k}\right) =  \gamma_k \mathbb{E}\left( \tilde{G}^S_{kk}(z)\right) + 
\tau_{i,k}$$
with $\vert \tau_{i,k}\vert \leq \frac{1}{p} \left( \vert z \vert +1 \right) 
P(\vert \Im z \vert^{-1}) $ and thus, using equation (\ref{eqinter}), there exists  a polynomial $Q$ with
nonnegative coefficients and a sequence of nonnegative numbers $a_N$ converging towards zero when $N$ goes to infinity such that for any $k=1,\ldots,N$ ,
$$\mathbb{E}\left( \tilde{G}^S_{kk}(z)\right)= \frac{1}{z} +  \gamma_k\mathbb{E}\left( \tilde{G}^S_{kk}(z)\right)
\mathbb{E}\left( \tr_p G_{\underline{M}_N^S}(z) \right) + \xi_k$$
with $$\vert \xi_k\vert \leq \left( \vert z \vert +1 \right)^2
Q(\vert \Im z \vert^{-1})a_N.$$
Thus \begin{equation}\label{eqinterm}
\left\{z-  \gamma_k z  \mathbb{E}\left( \tr_p G_{\underline{M}_N^S}(z) \right) \right\}
    \mathbb{E}\left( \tilde{G}^S_{kk}(z)\right)=1 + z \xi_k.
    \end{equation}
    Using  the resolvent identity $$ zG_{\underline{M}_N^S}(z)=I_p +  \underline{M}_N^SG_{\underline{M}_N^S}(z)$$
    we can easily see that 
    $$\Im \left[ z \tr_p G_{\underline{M}_N^S}(z) \right]=-(\Im z) \tr_p G_{\underline{M}_N^S}(z)^* \underline{M}_N^SG_{\underline{M}_N^S}(z).$$
    Hence \\
    
    $
    \vert \Im \left\{z-  \gamma_k z  \mathbb{E}\left( \tr_p G_{\underline{M}_N^S}(z) \right) \right\}\vert$
    \begin{equation} \label{image} = \vert \Im z \vert \left\{ 1+  \gamma_k\tr_p G_{\underline{M}_N^S}(z)^* \underline{M}_N^SG_{\underline{M}_N^S}(z) \right\} \geq \vert \Im z \vert.\end{equation}
    Thus (\ref{eqinterm}) yields that 
    $$\mathbb{E}\left( \tilde{G}_{kk}^S(z)\right)= \frac{1}{z-  \gamma_k z  \mathbb{E}\left( \tr_p G_{\underline{M}_N^S}(z)\right)} + \xi^{'}_k$$
  with   $$\vert \xi^{'}_k\vert \leq \frac{\vert z\vert }{\vert \Im z \vert}\left( \vert z \vert +1 \right)^2 
Q(\vert \Im z \vert^{-1})a_N.$$
Proposition \ref{estimresolvwish} readily follows since (see(\ref{baroupas})) we have
$$\mathbb{E}\left( \tr_p G_{\underline{M}_N^S}(z)\right)=\frac{N}{p} \mathbb{E}\left( \tr_N G_{{M}_N^S}(z)\right)
+ \frac{1-\frac{N}{p}}{z}.$$ $\Box$\\

To prove  Proposition \ref{estimdif} in the previous section, we need the following description, when the matrix $A_N$ and the measure $\nu$ satisfied Assumption A in the Introduction, of the  convergence of $g_{\mu_{A_N}}(z)$ towards $g_\nu(z)$
and of the convergence of ${\cal Z}_{\frac{N}{p},\mu_{A_N}}(z)$ towards ${\cal Z}(z)$ (dealing in the last case with measures on $[0;+\infty[$) where ${\cal Z}_{\frac{N}{p},\mu_{A_N}}(z)$ is defined by (\ref{defzcnu}) replacing $\nu$ by $\mu_{A_N}$ and $c$ by $\frac{N}{p}$. 
\begin{lemme}\label{cvAN} Under Assumption A,
there exists  polynomials $P_1$ and $P_2$ with nonnegative coefficients and  sequences $v_N(1)$ and $v_N(2)$ of positive numbers converging towards zero such that  for all $z \in \mathbb{C}\setminus \mathbb{R}$,
 \begin{equation}\label{hypA}\vert g_{\mu_{A_N}}(z) -g_\nu(z)\vert \leq P_1(\vert \Im z \vert^{-1}) v_N(1),\end{equation}
 \begin{equation}\label{hypA2}\vert {\cal Z}_{\frac{N}{p},\mu_{A_N}}(z) -{\cal Z}(z)\vert \leq (\vert z \vert +1)^2 P_2(\vert \Im z \vert^{-1}) v_N(2),\end{equation}
\end{lemme}
\noindent{\bf Proof}: Let us introduce $$\hat{\nu}_N= \frac{1}{N-r}\sum_{j=1}^{N-r} \delta_{\beta_j(N)}.$$
Let us fix $\epsilon >0$. According to the assumption (\ref{univconv}), for $N$ large all the $\beta_j(N)$
are in the set $\{x, d(x,\mbox{supp~} \nu)<{\epsilon} \}$. Moreover, $\{x, d(x,\mbox{supp~} \nu)< {\epsilon}\}$ may be covered by a finite number $n_\epsilon$ of disjoint intervals $I_i(\epsilon)$ with diameter smaller than $\epsilon$, of the form $]a_i(\epsilon);b_i (\epsilon)]$ where $a_i(\epsilon)$ and $b_i (\epsilon)$ are two continuity points of the distribution function of $\nu$.
Note that for any $i=1,\ldots, n_\epsilon,$
when $N$ goes to infinity,  $$\hat{\nu}_N (I_i(\epsilon)) \rightarrow \nu(I_i(\epsilon)).$$
Since $\vert \frac{1}{N}\sum_{i=1}^J \frac{1}{z-\theta_i}\vert \leq \frac{r}{N}\vert\Im z \vert^{-1}$,
and $\vert \left[\frac{1}{N}-\frac{1}{N-r}\right]\sum_{j=1}^{N-r} \frac{1}{z-\beta_j(N)}\vert \leq \frac{r}{N}\vert\Im z \vert^{-1}$,
we focus on the difference $g_{\hat{\nu}_N}(z) -g_\nu(z)$.
Similarly, since  $\vert \frac{1}{N}\sum_{i=1}^J \frac{\theta_i}{1-\theta_i z}\vert \leq \frac{r}{N}\vert\Im z \vert^{-1}$,
and $\vert \left[\frac{1}{N}-\frac{1}{N-r}\right]\sum_{j=1}^{N-r} \frac{\beta_j(N)}{1-\beta_j(N)z}\vert \leq \frac{r}{N}\vert\Im z \vert^{-1}$, we focus on the difference 
 ${\cal Z}_{\frac{N}{p},\hat{\nu}_N}(z) -{\cal Z}(z)$
where ${\cal Z}_{\frac{N}{p},\hat{\nu}_N}$  is defined by (\ref{defzcnu}) replacing $\nu$ by $\hat{\nu}_N$ and $c$ by $\frac{N}{p}$ .
\begin{eqnarray*}
g_{\hat{\nu}_N}(z) -g_\nu(z)&=& \sum_{i=1}^{n_\epsilon} \left\{ \frac{1}{N-r}\sum_{\beta_j(N) \in I_i(\epsilon)}\frac{1}{z-\beta_j(N)} - \int_{I_i(\epsilon)} \frac{1}{z-x} d\nu(x)\right\}\\
&=& \sum_{i,\nu(I_i(\epsilon))=0}  \frac{1}{N-r}\sum_{\beta_j(N) \in I_i(\epsilon)}\frac{1}{z-\beta_j(N)}\\
&+&\sum_{i,\nu(I_i(\epsilon))>0}  \frac{1}{N-r}\sum_{\beta_j(N) \in I_i(\epsilon)}\frac{1}{\nu(I_i(\epsilon))}
\int_{I_i(\epsilon)} \left( \frac{1}{z-\beta_j(N)} -  \frac{1}{z-x}\right) d\nu(x)\\
&+& \sum_{i,\nu(I_i(\epsilon))>0} \left( \frac{\hat{\nu}_N (I_i(\epsilon))} {\nu(I_i(\epsilon))} - 1 \right) 
\int_{I_i(\epsilon)}   \frac{1}{z-x}d\nu(x)\\
&=& \Delta_1+\Delta_2+\Delta_3.
\end{eqnarray*}
where
$$\vert \Delta_1\vert \leq \sum_{i,\nu(I_i(\epsilon))=0} \hat{\nu}_N (I_i(\epsilon)) \vert \Im z \vert^{-1},$$
$$\vert \Delta_2\vert \leq \epsilon \sum_{i,\nu(I_i(\epsilon))>0} \hat{\nu}_N (I_i(\epsilon)) \vert \Im z \vert^{-2}\leq \epsilon \vert \Im z \vert^{-2},$$
$$\vert \Delta_3\vert \leq \sum_{i,\nu(I_i(\epsilon))>0} \vert \hat{\nu}_N (I_i(\epsilon))- {\nu} (I_i(\epsilon)) \vert \vert \Im z \vert^{-1}.$$
Hence $$\vert g_{\hat{\nu}_N}(z) -g_\nu(z)\vert \leq \left(\vert \Im z \vert^{-2}+ \vert \Im z \vert^{-1}\right)\left(\epsilon +\sum_{i=1}^{n_\epsilon} \vert \hat{\nu}_N (I_i(\epsilon))- {\nu} (I_i(\epsilon)) \vert\right)$$
and then $$\limsup_{N\rightarrow + \infty} \sup_{z\in \mathbb{C}\setminus \mathbb{R}} \left\{\left(\vert \Im z \vert^{-2}+\vert \Im z \vert^{-1}\right)^{-1}
\vert g_{\hat{\nu}_N}(z) -g_\nu(z)\vert\right\} \leq \epsilon.$$
Since this is true for any $\epsilon >0$, we get that
$$\lim_{N\rightarrow + \infty} \sup_{z\in \mathbb{C}\setminus \mathbb{R}} \left\{\left(\vert \Im z \vert^{-2}+ \vert \Im z \vert^{-1}\right)^{-1}
\vert g_{\hat{\nu}_N}(z) -g_\nu(z)\vert\right\}=0$$
which yields (\ref{hypA}).\\
Now, since moreover $\vert {\cal Z}_{\frac{N}{p},\hat{\nu}_N}(z)- {\cal Z}_{c,\hat{\nu}_N}(z)\vert \leq \vert \frac{N}{p}-c\vert 
\vert \Im z \vert^{-1}$ we will study  ${\cal Z}_{c,\hat{\nu}_N}(z)-{\cal Z}(z)$.
Similarly,
\begin{eqnarray*}
\frac{1}{c}[{\cal Z}_{c,\hat{\nu}_N}(z) -{\cal Z}(z)]&=& \sum_{i=1}^{n_\epsilon} \left\{ \frac{1}{N-r}\sum_{\beta_j(N) \in I_i(\epsilon)}\frac{\beta_j(N)}{1-\beta_j(N)z} - \int_{I_i(\epsilon)} \frac{x}{1-xz} d\nu(x)\right\}\\
&=& \sum_{i,\nu(I_i(\epsilon))=0}  \frac{1}{N-r}\sum_{\beta_j(N) \in I_i(\epsilon)}\frac{\beta_j(N)}{1-\beta_j(N)z}\\
&+&\sum_{i,\nu(I_i(\epsilon))>0}  \frac{1}{N-r}\sum_{\beta_j(N) \in I_i(\epsilon)}\frac{1}{\nu(I_i(\epsilon))}
\int_{I_i(\epsilon)} \left( \frac{\beta_j(N)}{1-\beta_j(N)z} -  \frac{x}{1-xz}\right) d\nu(x)\\
&+& \sum_{i,\nu(I_i(\epsilon))>0} \left( \frac{\hat{\nu}_N (I_i(\epsilon))} {\nu(I_i(\epsilon))} - 1 \right) 
\int_{I_i(\epsilon)}   \frac{x}{1-xz}d\nu(x)\\
&=& \Delta_1+\Delta_2+\Delta_3.
\end{eqnarray*}
where
$$\vert \Delta_1\vert \leq \sum_{i,\nu(I_i(\epsilon))=0} \hat{\nu}_N (I_i(\epsilon)) \vert \Im z \vert^{-1},$$
\begin{eqnarray*}\vert \Delta_2\vert & \leq &\sum_{i,\nu(I_i(\epsilon))>0} \left\{ \frac{1}{N-r}\sum_{\beta_j(N) \in I_i(\epsilon)}\frac{1}{\nu(I_i(\epsilon))}
\int_{I_i(\epsilon)} \frac{\vert \beta_j(N)-x\vert }{\vert z\vert ^2(\vert \frac{1}{z}-\beta_j(N) \vert \vert \frac{1}{z}-x\vert  }  d\nu(x)\right\}\\&
\leq &
\epsilon \sum_{i,\nu(I_i(\epsilon))>0} \hat{\nu}_N (I_i(\epsilon)) \vert z\vert ^{-2} \left|\Im (\frac{1}{z}) \right|^{-2}\\
&\leq &\epsilon \vert z\vert ^2 \vert \Im z \vert^{-2},\end{eqnarray*}
$$\vert \Delta_3\vert \leq \sum_{i,\nu(I_i(\epsilon))>0} \vert \hat{\nu}_N (I_i(\epsilon))- {\nu} (I_i(\epsilon)) \vert \vert \Im z \vert^{-1}.$$
Hence $$\vert \frac{1}{c}[{\cal Z}_{c,\hat{\nu}_N}(z) -{\cal Z}(z)]\vert \leq \left( \vert z \vert^2 \vert \Im z \vert^{-2}+\vert \Im z \vert^{-1}\right)\left(\epsilon +\sum_{i=1}^{n_\epsilon} \vert \hat{\nu}_N (I_i(\epsilon))- {\nu} (I_i(\epsilon)) \vert\right)$$
and then $$\limsup_{N\rightarrow + \infty} \sup_{z\in \mathbb{C}\setminus \mathbb{R}} \left\{\left(\vert z \vert^2  \vert \Im z \vert^{-2}+\vert \Im z \vert^{-1}\right)^{-1}
\vert \frac{1}{c}[{\cal Z}_{c,\hat{\nu}_N}(z) -{\cal Z}(z)]\vert\right\} \leq \epsilon.$$
Since this is true for any $\epsilon >0$, we get that
$$\lim_{N\rightarrow + \infty} \sup_{z\in \mathbb{C}\setminus \mathbb{R}} \left\{\left( \vert z \vert^2 \vert \Im z \vert^{-2}+ \vert \Im z \vert^{-1}\right)^{-1}
\vert {\cal Z}_{c,\hat{\nu}_N}(z) -{\cal Z}(z) \vert\right\}=0$$
and (\ref{hypA2}) follows.
$\Box$\\

In the sample covariance matrix setting we will need the following upper bound of $\frac{1}{\Im (\underline{g}_N^S(z))}$
where $$\underline{g}_N^S(z)=\mathbb{E}\left( \tr_p G_{\underline{M}_N^S}(z)\right)=\frac{N}{p} \mathbb{E}\left( \tr_N G_{{M}_N^S}(z)\right)
+ \frac{1-\frac{N}{p}}{z},$$
with $$\underline{M}_N^S=\frac{1}{p} B_N^* A_N B_N.$$
\begin{lemme}\label{MajinversegN} There exists a constant $\underline{C}$ such that for any $z$ in $\mathbb{C}\setminus \mathbb{R}$,
$$\left|\frac{1}{\Im (\underline{g}_N^S(z))}\right| \leq \underline{C}(1+ \vert z \vert )^2 \vert \Im z \vert^{-1}$$
\end{lemme}
\noindent {\bf Proof}: Note that $$\vert \Im (\underline{g}_N^S(z)) \vert = \vert \Im z \vert \mathbb{E} \left[ \int_{\mathbb{R}} \frac{d \mu_{\underline{M}_N}(x)}{\vert z-x\vert^2} \right].$$
Now, for any $x$ in the spectrum of $\underline{M}_N^S$, 
$$\vert z-x\vert^2 \leq 2 (\vert z \vert^2 + \Vert \underline{M}_N^S \Vert^2) \leq 2 (\vert z \vert^2 + \Vert A_N \Vert^2
\Vert \frac{1}{p} B_N^* B_N \Vert^2 )$$
\noindent so that, with $C= \sup_N \Vert A_N \Vert$,
$$\mathbb{E} \left[ \int_{\mathbb{R}} \frac{d \mu_
{\underline{M}_N}(x)}{\vert z-x\vert^2} \right] \geq \mathbb{E} \left[ \frac{1}{ 2 (\vert z \vert^2 + C^2
\Vert \frac{1}{p} B_N^* B_N \Vert^2 )} \right].$$
According to Theorem 5.11 in \cite{BaiSil06}, $\left\| \frac{1}{p} B_N^* B_N \right\| =c(1+ \frac{1}{\sqrt{c}})^2+o_{a.s,N}(1)$ so that by the dominated convergence Theorem we can deduce that for all large $N$,
$$\mathbb{E} \left[ \int_{\mathbb{R}} \frac{d \mu_
{\underline{M}_N}(x)}{\vert z-x\vert^2} \right] \geq  \frac{1}{ 2 (\vert z \vert^2 + C_1
 )} $$ \noindent where $C_1 > [C c(1+ \frac{1}{\sqrt{c}})^2]^2.$ Therefore 
 $$\left|\frac{1}{\Im (\underline{g}_N^S(z))}\right| \leq \frac{2 (\vert z \vert^2 + C_1
 )}{\vert \Im z\vert }.\,  $$
 so that Lemma \ref{MajinversegN} readily follows. $\Box$

\section{Appendix}
\subsection{Poincar\'e inequality and concentration inequalities}
We first derive in this section concentration inequalities based on the Poincar\'e inequality.
We refer the reader to the book \cite{Tou}. 
A probability measure $\mu$ on $\mathbb{R}$ is said to satisfy the Poincar\' e inequality with constant $C_{PI}$ if
 for any
${\cal C}^1$ function $f: \R\rightarrow \C$  such that $f$ and
$f' $ are in $L^2(\mu)$,
$$\mathbf{V}(f)\leq C_{PI}\int  \vert f' \vert^2 d\mu ,$$
\noindent with $\mathbf{V}(f) = \int \vert
f-\int f d\mu \vert^2 d\mu$. \\
 We refer the reader to \cite{BobGot99} for a characterization
 of  the measures on  $\mathbb{R}$ which satisfy a Poincar\'e inequality. 

\begin{remarque}\label{multiple} If the law of a random variable $X$ satisfies the Poincar\'e inequality with constant $C_{PI}$ then, for any fixed $\alpha \neq 0$, the law of $\alpha X$ satisfies the Poincar\'e inequality with constant $\alpha^2 C_{PI}$.\\
If a probability measure $\mu$ on $\mathbb{R}$ satisfies the Poincar\'e inequality with constant $C_{PI}$ then the product measure $\mu^{\otimes M}$ on $\mathbb{R}^M$ satisfies the Poincar\'e inequality with constant $C_{PI}$ in the sense that for any differentiable function $F$ such that $F$ and its gradient $\nabla F$ are in $L^2(\mu^{\otimes M})$,
$$\mathbf{V}(f)\leq C_{PI} \int \Vert \nabla F \Vert_2 ^2 d\mu^{\otimes M}$$
\noindent with $\mathbf{V}(f) = \int \vert
f-\int f d\mu^{\otimes M} \vert^2 d\mu^{\otimes M}$ (see Theorem 2.5 in  \cite{GuZe03}) .
\end{remarque}

An important consequence of the Poincar\'e inequality is the following concentration result.
\begin{lemme}\label{Herbst}{ Lemma 4.4.3 and Exercise 4.4.5 in \cite{AGZ09} or Chapter 3 in  \cite{Ledoux01}. }
Let $\mathbb{P}$ be a probability measure on $\mathbb{R^M}$ which satisfies a Poincar\'e inequality with constant $C_{PI}$. Then there exists $K_1>0$ and $K_2>0$ such that,  for any  Lipschitz function $F$  on $\mathbb{R}^M$ with Lipschitz constant $\vert F \vert_{Lip}$,
$$\forall \epsilon> 0,  \, \mathbb{P}\left( \vert F-\E_{\mathbb{P}}(F) \vert > \epsilon \right) \leq K_1 \exp\left(-\frac{\epsilon}{K_2 \sqrt{C_{PI}} \vert F \vert_{Lip}}\right).$$
\end{lemme}

\subsection{Technical tools}
 We need the following result on the extension of Lipschitz functions on $\mathbb{R}$ to the Hermitian matrices.
\begin{lemme}\label{extlipschitz}(see \cite{D})
Let $f$ be a real $C_{\cal L}$-Lipschitz function on $\mathbb{R}$. Then its extension on the $N\times N$
Hermitian matrices is $C_{\cal L}$-Lipschitz with respect to the norm $\Vert M\Vert_2=\{Tr (MM^*)\}^{\frac{1}{2}}$.
\end{lemme}
\noindent {\bf Proof:} Let $A$ and $B$ be $N\times N$
Hermitian matrices. Let us consider their  spectral decompositions 
$$A=\sum_i \lambda_i(A)P^{(A)}_i$$ \noindent and 
$$B=\sum_i \lambda_i(B)P^{(B)}_i.$$ \noindent We have
\begin{eqnarray*}
\Vert f(B) -f(A) \Vert^2_2 &=& Tr \left( \sum_i f(\lambda_i(A))P^{(A)}_i -\sum_i f(\lambda_i(B))P^{(B)}_i \right)^2\\
&=&Tr \left( \sum_{i} f(\lambda_i(A))^2P^{(A)}_i +\sum_{j} f(\lambda_j(B))^2P^{(B)}_j \right)
\\ & &-2 \sum_{i,j} f(\lambda_i(A))f(\lambda_j(B))Tr (P^{(A)}_i P^{(B)}_j)\\
 &=&Tr \left( \sum_{ij} (f(\lambda_i(A))^2P^{(A)}_i P^{(B)}_j+\sum_{i,j} (f(\lambda_j(B))^2P^{(A)}_i P^{(B)}_j \right)
\\ & & -2 \sum_{i,j} f(\lambda_i(A))f(\lambda_j(B))Tr (P^{(A)}_i P^{(B)}_j)\\
&=& \sum_{i,j} (f(\lambda_i(A))-f(\lambda_j(B))^2Tr( P^{(A)}_i P^{(B)}_j).
\end{eqnarray*}
Now, since $Tr( P^{(A)}_i P^{(B)}_j)\geq 0$, we can deduce that 
$$\Vert f(B) -f(A) \Vert_2^2 \leq \sum_{i,j} C_{\cal L}^2(\lambda_i(A)-\lambda_j(B))^2Tr( P^{(A)}_i P^{(B)}_j)= C_{\cal L}^2 \Vert B-A \Vert_2^2. \Box$$

We recall here some useful properties of the resolvent (see \cite{KKP96,CD07}).
\begin{lemme} \label{lem0}
For a $N \times N$ Hermitian or symmetric matrix $M$, 
for any $z \in \C\setminus {\rm Spect}(M)$, 
we denote by $G(z) := (zI_N-M)^{-1}$ the resolvent of $M$.\\
Let $z \in \C\setminus \R$, 
\begin{itemize}
\item[(i)] $\Vert G(z) \Vert \leq |\Im z|^{-1}$ where $\Vert . \Vert$ denotes the operator norm. 
\item[(ii)] $\vert G(z)_{ij} \vert \leq |\Im z|^{-1}$ for all $i,j = 1, \ldots , N$. 
\item[(iii)] Let $z \in \C$ such that $|z| > \Vert M \Vert$; we have
$$\Vert G(z) \Vert \leq \frac{1}{|z| - \Vert M \Vert}.$$
\end{itemize}
\end{lemme}
We recall here the following classical result due to
Weyl.

\begin{lemme}{(cf. Theorem 4.3.7 of \cite{HJ})} \label{Weyl}
Let B and C be two $N \times N$ Hermitian matrices. For any pair of integers
$j,k$ such that $1 \leq j,k\leq N$ and $j+k \leq N+1$, we have
$$\lambda_{j+k-1} (B+C) \leq \lambda_{j}(B) + \lambda_{k}(C).$$
For any pair of integers
$j,k$ such that $1 \leq j,k\leq N$ and $j+k \geq N+1$, we have
$$ \lambda_j(B) + \lambda_k(C) \leq \lambda_{j+k-N} (B+C).$$
\end{lemme}

The following  result on quadratic forms  is of basic use in the sample covariance matrix setting. Note that, a complex random variable $x$ will be said {\it standardized} if $\mathbb{E}(x)=0$ and $\mathbb{E}(\vert x \vert^2)=1$.
\begin{proposition}\label{QF}(Lemma 2.7  \cite{BaiSil98})
Let $B=(b_{ij})$ be a $N \times N$  matrix and
$Y_N$ be a vector of size $N$ which contains i.i.d
standardized entries with bounded fourth
moment. Then there is a constant $K>0$ such that $$\mathbb E\vert  Y_N ^* B Y_N - {\rm{Tr}} B\vert^2 \leq K \mathbb \Tr (BB^*).$$
\end{proposition}

The following technical lemma is fundamental  in this paper.
We refer the reader to the Appendix of \cite{CD07} 
where it is proved using the ideas of \cite{HaaThor05}.
\begin{lemme}\label{HT} Let $h$ be an analytic function on $\C\setminus \R$ which satisfies
\begin{equation*}\label{nestimgdif}
\vert h(z)\vert \leq (\vert z\vert +K)^\alpha P(\vert \Im z\vert ^{-1})
\end{equation*} and $\varphi$ be in $\cal C^\infty (\R, \R)$ with compact support. Then,
$$\limsup _{y\rightarrow 0^+}\vert \int _\R\varphi (x)h(x+iy)dx\vert < + \infty.$$
\end{lemme}

\noindent {\bf Acknowledgments:}
I am grateful to Charles Bordenave 
for  useful discussions.  I would like to thank the anonymous referees  for their careful reading and their pertinent comments which led to an overall improvement of the paper.

 
\def\cprime{$'$}

\end{document}